\documentclass[12pt,reqno]{amsart}
\usepackage{amssymb}
\usepackage{hyperref}
\usepackage{amsmath}
\usepackage{amsthm}
\usepackage{geometry}
\usepackage{amsfonts, adjustbox, amssymb, amscd}
\numberwithin{equation}{section}

\usepackage{bm}
\usepackage{verbatim}
\usepackage{mathrsfs}
\usepackage{graphicx}
\usepackage{tikz-cd}
\usepackage{comment}
\usepackage{tikz}
\usetikzlibrary{shapes.geometric, arrows}
\tikzstyle{startstop} = [rectangle, rounded corners, minimum width=3cm, minimum height=1cm,text centered, draw=black]
\tikzstyle{io} = [trapezium, trapezium left angle=70, trapezium right angle=110, minimum width=3cm, minimum height=1cm, text centered, text width=7cm, draw=black]
\tikzstyle{process} = [rectangle,  minimum height=1cm, text centered, text width=7cm, draw=black]
\tikzstyle{decision} = [diamond, minimum width=3cm, minimum height=0.5cm, text centered, draw=black]
\tikzstyle{arrow} = [thick,->,>=stealth]

\usepackage{subcaption}
\usepackage{listings}
\usepackage{subfiles}
\usepackage[toc,page]{appendix}
\usepackage{mathtools}
\usepackage{comment}
\usepackage{enumerate}
\usepackage{enumitem}
\usepackage[all]{xy}

\usepackage{graphicx}
\usepackage{appendix}
\usepackage{hyperref}
\hypersetup{
    colorlinks=true,
    citecolor=red,
    linkcolor=blue,
    filecolor=magenta,      
    urlcolor=red,
}
\usepackage[dvipsnames]{xcolor}

\newcommand{\Spec}{\mathrm{Spec}}

\newcommand{\Cc}{\mathbb{C}}

\newcommand{\Qq}{\mathbb{Q}}

\newcommand{\ga}{\mathbb{G}_a}
\newcommand{\gm}{\mathbb{G}_m}

\newcommand{\Bir}{\operatorname{Bir}}
\newcommand{\Bim}{\operatorname{Bim}}

\newcommand{\Aut}{\operatorname{Aut}}

\newcommand{\lct}{\operatorname{lct}}

\newcommand{\proj}{\operatorname{Proj}}
\newcommand{\spec}{\operatorname{Spec}}

\newcommand{\Supp}{\operatorname{Supp}}

\newcommand{\Diff}{\operatorname{Diff}}

\newcommand{\mult}{\operatorname{mult}}

\newcommand{\lf}{\lfloor}
\newcommand{\rf}{\rfloor}

\newcommand{\Oo}{\mathcal{O}}


\newcounter{parentnumber}

\newtheorem{thm}{Theorem}[section]

\newtheorem{lem}[thm]{Lemma}
\newtheorem{lemma}[thm]{Lemma}

\newtheorem{proposition}[thm]{Proposition}
\newtheorem{ex}[thm]{Example}
\newtheorem{claim}{Claim}

\newtheorem*{lem*}{Lemma}
\newtheorem*{prop*}{Proposition}

\newtheorem{innercustomthm}{Theorem}

\theoremstyle{definition}

\newtheorem{rmk}[thm]{Remark}
\newtheorem{definition}[thm]{Definition}

\title{Abundance for semi log canonical threefolds in positive characteristic}
\author{Sheng Qin}
\address{School of Mathematical Sciences, Peking University, Beijing, China.}
\email{qinsheng@stu.pku.edu.cn}
\date{\today}

\begin{document}

\begin{abstract}
We work over an algebraically closed field $k$ of characteristic $p> 42$. We prove the abundance conjecture for semi log canonical threefolds over $k$. As a key step, we prove the finiteness of $B$-pluricanonical representations for surfaces over $k$.
\end{abstract}

\maketitle
\tableofcontents

\section{Introduction}

The abundance conjecture is a very important problem in the birational classification of algebraic varieties.

\medskip
\noindent\textbf{Abundance conjecture.} Let $(X, B)$ be a projective semi log canonical pair of dimension $3$ over a perfect field $k$ of characteristic $>0$. If $K_X + B$ is nef, then it is semi-ample.

\medskip
In this paper, we prove the abundance conjecture for semi log canonical threefolds in characteristic $p>42$. This is a generalization of the abundance conjecture for log canonical threefolds proved by Bernasconi, Witaszek and Xu in their work \cite{bernasconi2026abundancethreefoldspositivecharacteristic}.

\begin{thm}[Theorem \ref{cor:slc-threefold-abundance}]\label{main thm abundance}
    Let $(X, B)$ be a projective semi log canonical pair of dimension $3$ over an algebraically closed field $k$ of characteristic $p>42$. If $K_X + B$ is nef, then it is semi-ample.
\end{thm}


The notion of semi log canonical singularity was first introduced in \cite{kollar1988threefolds} in order to compactify the moduli of surfaces and has been applied in the moduli of varieties in higher dimensions (see Koll\'ar's book \cite{kollar-moduli-hd}). There are other reasons to consider semi log canonical singularities. For example, given a divisorial log terminal pair $(X,B)$, the reduced part of the boundary appears to be only semi-normal instead of normal, and the adjunction $(\lfloor B\rfloor,\Diff(B-\lfloor B\rfloor))$ is proved to be a semi divisorial log terminal pair in characteristic zero. Also, the abundance for semi log canonical $(n-1)$-folds is viewed as an important step in proving the abundance for log canonical $n$-folds (see \cite[Section 7]{abundanceforthreefolds}).

In characteristic zero, the abundance for semi log canonical threefolds was proved in \cite{Fuj00},  where the author first explicitly introduced $B$-birational maps (see \cite[Definition 1.5]{Fuj00}) and used the finiteness of $B$-pluricanonical representations to descend semi-ampleness along the normalization. Tanaka applied this approach to prove abundance for semi log canonical surfaces over an algebraically closed field of positive characteristic in \cite{tanaka2015abundancetheoremsemilog}.

There is another approach to the abundance for semi log canonical varieties.  In \cite{hacon2012finitenessbrepresentationsemilogcanonical}, the authors use the theory of sources and springs developed by Koll\'ar (see \cite{kollár2012sourceslogcanonicalcenters}), which allows one to construct a quotient space with set-theoretic equivalence relation. The finiteness of $B$-pluricanonical representations is used to prove the finiteness of a set-theoretic equivalence relation. In \cite{Posva_2023_slc}, the author uses similar strategies to prove the abundance conjecture for semi log canonical surfaces over arbitrary fields of positive characteristic.

These two approaches are actually closely related, and the finiteness of $B$-pluricanonical representations is crucial in both approaches. Let us briefly review the history of this finiteness result. We first recall the following classical result. For a compact complex manifold $X$, let $\Bim(X)$ denote the group of bimeromorphic self-maps of $X$.  We have a pluricanonical representation
$$\rho_m: \Bim(X)\to \operatorname{Aut}(H^0(X,mK_X))$$
for any positive integer $m$. Nakamura and Ueno first proved the following finiteness theorem for Moishezon manifolds under the assumption that the image of the pluricanonical mapping $\Phi_{mK_X}(X)$ is not a ruled variety (see \cite[Theorem 1(3)]{Nakamura_Ueno[1]}). Deligne then proved the theorem in full generality.

\begin{thm}[cf.{\cite[Theorem 14.10]{Ueno1975}}]
    Let $X$ be a compact complex manifold. Assume that $X$ is Moishezon. Then the group $\rho_m(\Bim(X))$ is finite.
\end{thm}

Fujino and Gongyo extended this theorem to the following setting for projective log canonical pairs with semi-ample log canonical divisor. Hacon and Xu gave another proof using Hodge theory and local systems.

\begin{thm}[cf.{\cite[Theorem 1.1]{Fujino_Gongyo_2014}, \cite[Theorem 1.2]{hacon2012finitenessbrepresentationsemilogcanonical}}]
    Let $(X,B)$ be a projective log canonical pair over the field $\Cc$ of complex numbers. Suppose that $m(K_X+B)$ is Cartier and that $K_X+B$ is semi-ample. Then the image of the $m$-th $B$-pluricanonical representation $\rho_m(\Bir(X,B))$ is finite.
\end{thm}
The classical complex proofs use analytic techniques, such as integration of differential forms or Hodge theory. Unfortunately, these do not have direct analogues in positive characteristic. So, little is known about the finiteness of pluricanonical representation in characteristic $p>0$. Moreover, the canonical bundle formula and the behavior of singularities are subtler in characteristic $p>0$. In this paper, we prove the following finiteness results for surfaces.


\begin{thm}[Theorem \ref{main thm: finiteness of pluri for peff lc surf & Cartier ind}]\label{main thm finiteness}
    Let $(S,B)$ be a projective log canonical surface over an algebraically closed field $k$ of characteristic $p >42$. Assume that $K_S+B$ is nef. Then the image of the $B$-pluricanonical representation
    \[
        \rho_{m}:\Bir(S,B)\longrightarrow
        \Aut\big( H^0(S,m(K_S+B))\big)
    \]
is finite for every positive integer $m$ such that $m(K_S+B)$ is Cartier.
\end{thm}

\medskip
\noindent\textbf{Outline of the proof of Theorem~\ref{main thm abundance}.}
We use Fujino's strategy in \cite{Fuj00} to reduce to Theorem \ref{main thm finiteness}. To prove this, the first step is to treat Kawamata log terminal pairs. Here, we explain the strategy to prove Theorem \ref{main thm: finite pluri repn for Kawamata log terminal}. Note that we do not assume nefness of $K_S+B$ in this situation (see also \cite[Theorem 3.4]{Fuj00}). By taking a terminal model, we reduce to the finiteness of $A$-pluricanonical representations (see Definition \ref{def: B-bir Bir()}). After further modifications, we may assume the coefficients of $B$ to be standard; that is,
$$B=\sum_i(1-1/n_i)B_i$$
for some irreducible divisors $B_i$.

We may assume $K_S+B$ is pseudo-effective. By Lemma \ref{emmp}, we may run a $(K_S+B)$-MMP and further assume that $K_S+B$ is nef. By the abundance for surfaces (see \cite[Theorem 1.2]{abundanceforlcsurf}), $K_S+B$ is semi-ample. Consider the induced fibration $f:S\to C$. Now we divide the proof into three cases according to $\kappa(K_S+B)$. In the case of $\kappa(K_S+B)=2$, we conclude by finiteness of the group $\Aut(S,B)$ (see \cite[Corollary 9.5]{Pat17}). In the case of $\kappa(K_S+B)=0$, the bound of global indices for two-dimensional log Calabi--Yau pairs proved in \cite{BWX26indices} allows us, under the assumption of $\operatorname{char} k>42$, to take a tame index-one-cover. Therefore, after taking the minimal resolution, we may reduce to the case $(S,B)=(S,0)$ with $S$ being a $K3$ surface or an abelian surface (see Lemma \ref{B bir lifted to cyclic cover}). Now we may conclude by Theorem \ref{thm:abelian-variety-pluricanonical-finite} and \ref{main thm for smooth K3 surface}.

In the case of $\kappa(K_S+B)=1$, $C$ is a projective smooth curve. Note that every element in $\Aut(S,B)$ acts on the fibration equivariantly. This induces a homomorphism
$$l:\Aut(S,B)\to\Aut(C).$$
We may first reduce the finiteness of $A$-pluricanonical representations for $(S,B)$ to the finiteness of the image of $l$ (see the paragraphs above Claim \ref{claim:image} and those above \cite[Theorem 3.1]{hacon2012finitenessbrepresentationsemilogcanonical}). Suppose, for contradiction, that the image of $l$ is infinite. Then we investigate the canonical bundle formula. We may further assume that $f$ is a $\mathbb{P}^1$-bundle or a minimal elliptic fibration.

If $f$ is a $\mathbb{P}^1$-bundle, we apply the canonical bundle formula proved by \cite[6.2]{cascini2015base} (see also \cite[Theorem 4.8]{DH16}, \cite[Lemma 3.1]{Witaszek2021}) to show that the moduli part is Q-linearly trivial (see Lemma \ref{canonical bundle formula with trivial moduli}) and then reduce the question to the one-dimensional base $C$. If $f$ is a minimal elliptic fibration, we apply the classification results in \cite{Gu2020SurfaceFW}. Computation of the canonical bundle $K_S$ shows that $K_S+B\leq 0$, which leads to a contradiction.

\medskip
\noindent\textbf{Acknowledgements.}
The author would like to express his extreme gratitude to Zheng Xu for suggesting this problem and for his generous and ongoing academic guidance, encouragement and support. The author also thanks him for his helpful comments and suggestions on the manuscript.

The author thanks Haofeng Zhang for his helpful comments on an earlier draft of this paper. The author is grateful to Yuan Lu for helpful suggestions on the formatting of the paper, and thanks Caucher Birkar, Ruicheng Hu, Xintong Jiang, Jakub Witaszek and Tianle Yang for useful discussions. The author also thanks Chenyang Xu for useful discussions on a finite base change technique the author first learned in \cite{hacon2013birational,hacon2018boundedness}.

The author would also like to thank the organizers of the Algebra and Number Theory Summer School in 2025, after which this project began, for creating a stimulating academic environment and providing this valuable research opportunity.

\medskip
\noindent\textbf{AI disclosure.}
This project began after the summer of 2025. The author completed the main body of the proofs during the autumn semester of 2025. At this point, ChatGPT was used to supply Theorem \ref{thm:abelian-variety-pluricanonical-finite}. Also, for Lemma \ref{emmp}, ChatGPT suggested using Artin's contraction theorem to contract the negative part of a Zariski decomposition, obtaining a model with semi-ample log canonical divisor. The construction yields an algebraic space. Building on this suggestion, the author developed an argument using MMP to obtain the required projective variety and completed the proof of Lemma \ref{emmp}. ChatGPT was also used to search the literature, including the paper \cite{Gu2020SurfaceFW}. The subsequent work described below concerns revising, correcting these arguments later.

Shortly before the summer of 2026, the author used the newest ChatGPT model at that time for review and it located an example in \cite{katsura1984elliptic} that ruled out the proposed general argument for Theorem \ref{main thm of kappa one} using the canonical bundle formula. The author then decided to treat elliptic fibrations and $\mathbb P^1$-bundles separately and finished the argument in July 2026.

After the summer of 2026, the author submitted the draft to ChatGPT 6.1 for review. It produced a counterexample showing that the original boundary-modification argument for the $\mathbb P^1$-bundle case did not work. The author then guessed that the canonical-bundle-formula argument might still be valid under the particular hypotheses of this case, and used ChatGPT 6.1 to locate the relevant literature.

While preparing the final manuscript, the author used ChatGPT 6.2 to proofread and refine the exposition of the proofs. The author takes full responsibility for all mathematical claims, proofs, references, and the final text.

\section{Preliminaries}
In this section, we recall some basic notation and results. We always work over an algebraically closed field $k$ of characteristic $p>0$.

\begin{itemize}
    \item A \emph{variety} is a separated, reduced, and equi-dimensional scheme of finite type over $k$. Note that a variety in our sense might be reducible.
    \item A \emph{$\Qq$-divisor} $D = \sum d_jD_j$ on a scheme $X$ is a finite formal sum such that $d_j\in\Qq$ and $D_j$ is a prime divisor for each $j$ with $D_i\neq D_j$ for $i\neq j$. We define the \emph{round-down} of $D$ to be $\lfloor D\rfloor := \sum_j \lfloor d_j\rfloor D_j$.
    \item We say that two $\Qq$-divisors $D_1,D_2$ are \emph{$\Qq$-linearly equivalent}, denoted by $D_1\sim_\Qq D_2$, if there exists a positive integer $N$, such that $ND_1,ND_2$ are both integral divisors such that $ND_1$ and $ND_2$ are linearly equivalent.
    \item If $X$ is normal and $D$ is an integral Weil divisor on $X$, we may write $H^0(X,D)$ to denote $H^0(X,\mathcal O_X(D))$.
    \item Let $D$ be a $\Qq$-divisor on a normal projective variety $X$. The \emph{Kodaira dimension} $\kappa(D)$ is defined to be the smallest non-negative integer $\kappa$ such that
        $$\dim H^0(X,\lfloor mD\rfloor)\leq cm^\kappa$$
for some positive constant $c$ and for all large positive integers $m$ if the left-hand side is not always zero, and is equal to $-\infty$ otherwise. It is known that if $D$ is semi-ample, then $\kappa(D)$ equals the dimension of the image of $X$ under the morphism defined by the linear system $|mD|$ for a sufficiently divisible positive integer $m$.
\end{itemize}

\subsection{Pairs and singularities}
In this subsection, we recall some basic definitions and notation in the theory of the minimal model program. We refer to \cite{Kollár_Mori_1998} for a more detailed treatment.

\begin{definition}\label{sing of pairs}
    A \emph{sub-pair} $(X, B)$ consists of a normal variety $X$ and a $\Qq$-divisor $B$ such that $K_X + B$ is $\Qq$-Cartier. If $B \geq 0$, then $(X, B)$ is called a \emph{pair}. If all coefficients of $B$ lie in $[0,1]$, then $B$ is called a \emph{boundary} of $X$. We say that a divisor $B$ has \emph{standard coefficients} if the coefficients of $B$ belong to the set $\{1-\frac{1}{n}|n\in\mathbb{N}_{>0}\}$.
    
     A pair $(X, B)$ is called \emph{log smooth} if $X$ is smooth and the divisor $B$ has simple normal crossing support.
     
    A projective birational morphism $g : Y \rightarrow X$ is called a \emph{log resolution} of a pair $(X, B)$ if $Y$ is smooth and the strict transform of $B$ union the exceptional locus is a divisor with simple normal crossing support.
    
    Let $f: Y \to X$ be any log resolution of $(X,B)$, and let
\[
K_Y + B_Y:= f^*(K_X + B).
\]
The \emph{log discrepancy} of a prime divisor $D$ on $Y$ with respect to $(X, B)$ is $1-\mult_D B_Y$ and it is denoted by $a_D( X, B)$. We say that $(X, B)$ is \emph{Kawamata log terminal} (resp. \emph{log canonical}) if $a_D( X, B) > 0$, (resp. $\geq 0$) for every log resolution $ f : Y\to X$ and every prime divisor $D$ on $Y$. We say that $(X, B)$ is \emph{terminal} if $\lfloor B\rfloor=0$ and $a_E( X, B)>1$ for every log resolution $f : Y \to X$ and every exceptional prime divisor $E$ on $Y$. We say that $(X, B)$ is \emph{divisorial log terminal} if it is log canonical and $a_D( X, B) > 0$ for any exceptional prime divisor $D \subset Y$ over $X$ for some log resolution $f : Y \to X$.

\begin{definition}
    Let $(X, B) $ be a sub-pair and let $E$ be a divisor over $X$. We say that $E$ is a \emph{non-klt place} for the sub-pair if $a_E(X, B) \le 0$; in this case, we say that the center of the divisor $E$ in $X$, $ c_X(E) \subset X$, is a \emph{non-klt center} for the sub-pair. The \emph{non-klt locus} is defined as the union of all the non-klt centers of $(X, B)$ and we denote it by $\operatorname{Nklt}(X, B)$.
\end{definition}
\end{definition}

\subsection{Pluricanonical representation}
We introduce the notion of \emph{$B$-pluricanonical representations} in this subsection.
\begin{definition}\label{def: B-bir Bir()}
   Let $(X,B)=\coprod_{i=1}^n(X_i,B_i)$, $(Y,B')=\coprod_{i=1}^n(Y_i,B'_i)$ be two normal projective sub-pairs with their decompositions into irreducible components. We say that a proper birational map $f : (X, B)\dashrightarrow(Y, B')$ is \emph{B-birational} if there exists a normal variety $W$ and proper birational morphisms $r_1:W\rightarrow X$ and $r_2:W\rightarrow Y$ such that we have the following commutative diagram
    \begin{center}
        \begin{tikzcd}
                          & W \arrow[ld, "r_1"'] \arrow[rd, "r_2"] &   \\
X \arrow[rr, "f", dashed] &                                             & Y
\end{tikzcd}
    \end{center}
    and
    $$r_1^*(K_X + B) = r_2^*(K_Y + B').$$
    We define
    $$\Bir(X,B):=\{f:(X,B)\dashrightarrow(X,B)| f\text{ is a $B$-birational map}\},$$
    $$\Aut(X,B):=\{f:X\rightarrow X|f \text{ is an automorphism and }f^*B=B\}.$$
    
For any positive integer $m$ such that $m(K_X+B)$ is Cartier, a $B$-birational map $f\in\Bir(X,B)$ acts on $H^0(X,m(K_X+B))$ as follows
$$ \alpha\mapsto f_*\alpha:={r_2}_*\circ r_1^*(\alpha)$$
for any section $\alpha\in H^0(X,m(K_X+B))$.

We define the induced representation of $\Bir(X,B)$ to be the \emph{$m$-th $B$-pluricanonical representation}, denoted by
    $$\rho_m:\Bir(X,B)\rightarrow \Aut(H^0(X,m(K_X+B))),$$
and the \emph{$m$-th $A$-pluricanonical representation} is defined to be the restricted representation of $\rho_m$ on $\Aut(X,B)$. For simplicity of notation, we still denote it by
    $$\rho_m:\Aut(X,B)\rightarrow \Aut(H^0(X,m(K_X+B))).$$
\end{definition}

\begin{ex}
(Quadratic transformation). Let $X = \mathbb{P}^2$ and let $B$ be the union of three general lines on $\mathbb{P}^2$. Let $\alpha \colon W \to X$ be the blow-up at the three intersection points of $B$ and let $\beta \colon W \to X$ be the blow-down of the strict transform of $B$ on $W$. Then we obtain the \textit{quadratic transformation} $\varphi:X\dashrightarrow X$
\[
\begin{tikzcd}
& W \arrow[dl, "\alpha"'] \arrow[dr, "\beta"] & \\
X \arrow[rr, dashed, "\varphi"'] & & X
\end{tikzcd}
\]
It follows that
\[
\alpha^*(K_X + B) = K_W+\alpha_*^{-1}B+E_1+E_2+E_3 =\beta^*(K_X + B),
\]
where $E_1, E_2, E_3$ are the exceptional divisors of the morphism $\alpha$. Therefore, $\varphi$ is a $B$-birational map of the pair $(X, B)$.
\end{ex}

The next lemma is used in Proposition \ref{reduced bdry of dlt 3 folds}.

\begin{lem}\label{decsend CY under MMP}
    Let \((X, B)\) and \((Z, B')\) be normal varieties with \(\mathbb{Q}\)-divisors such that \(K_X + B\) and \(K_Z + B'\) are \(\mathbb{Q}\)-Cartier. Let \(f : X \rightarrow R\) and \(h : Z \rightarrow R\) be proper surjective morphisms onto a normal variety \(R\) and \(p : X \dashrightarrow Z\) a birational map such that \(f = h \circ p\). Assume that
\begin{enumerate}
    \item \(p^{-1}\) has no exceptional divisors;
    \item \(B' = p_*B\);
    \item \(K_X + B \equiv_f 0\), \(K_Z + B' \equiv_h 0\).
\end{enumerate}
Then \(p\) is a B-birational map.
\end{lem}
\begin{proof}
    Let \(\beta : W \rightarrow X\), $\alpha = p \circ \beta : W \rightarrow Z$ be a common resolution. Let \(m > 0\) be a sufficiently divisible integer. We have linear equivalences
\[
-mK_W = -\alpha^*(m(K_Z + B')) - F,
\]
\[
mK_W = \beta^*(m(K_X + B)) + E.
\]
Adding the two we obtain
\[
\beta^*(m(K_X + B)) - \alpha^*(m(K_Z + B')) = F - E.
\]
By assumption, \(F - E\) is \(\alpha\)-exceptional and numerically \(\alpha\)-trivial. Thus \(F = E\) and $p$ is a $B$-birational map.
\end{proof}

\begin{rmk}\label{Main Thm equi under B-bir}
    Assume that $f:(X,B)\dashrightarrow(X',B')$ is $B$-birational, then the finiteness of $B$-pluricanonical representations for $(X,B)$ holds if and only if that holds for $(X',B')$.

    To see this, we take any $B$-birational map $g\in\Bir(X,B)$, and consider the composition $g':=f\circ g\circ f^{-1}\in\Bir(X')$. After taking a sufficiently high log resolution $W$, we have
    $$r_1'^*(K_{X'}+B')=r_1^*(K_{X}+B)=r_2^*(K_{X}+B)=r_2'^*(K_{X'}+B'),$$
thus $g'\in\Bir(X',B')$.
    \begin{center}

\begin{tikzcd}
                                                        & W \arrow[ldd, "r_1"'] \arrow[lddd, "r_2" description] \arrow[rdd, "r_1'"] \arrow[rddd, "r_2'" description] &                                   \\
                                                        &                                                                                                            &                                   \\
{(X,B)} \arrow[rr, "f", dashed] \arrow[d, "g"', dashed] &                                                                                                            & {(X',B')} \arrow[d, "g'", dashed] \\
{(X,B)} \arrow[rr, "f", dashed]                         &                                                                                                            & {(X',B')}                        
\end{tikzcd}
    \end{center}
For any section $\alpha\in H^0(X,m(K_X+B))$, 
    $$f_*(g_*\alpha)=f_*({r_2}_*\circ r_1^*(\alpha))=f_*({r_2}_*\circ r_1'^*(f_*\alpha))={r_2'}_*\circ r_1'^*(f_*\alpha)=g'_*(f_*\alpha)$$
    
Given the isomorphism between $H^0(X,m(K_X+B))$ and $H^0(X',m(K_{X'}+B'))$ induced by $f_*$ and $f^{-1}_*$, and the isomorphism between $\Bir(X,B)$ and $\Bir(X',B')$ induced by conjugation, we obtain that the finiteness of $B$-pluricanonical representations for one of these two pairs implies that of the other one.
\end{rmk}

We prove the following finiteness of $B$-pluricanonical representations for curves.

\begin{thm}\label{1-dim subklt fin rep}
    Let $(C,B)$ be a projective one-dimensional log canonical sub-pair. Assume $\deg(K_C+B)\geq0$. Then the image of the $B$-pluricanonical representation
    \[
        \rho_{m}:\Aut(C,B)\longrightarrow
        \Aut\big( H^0(C,m(K_C+B))\big)
    \]
is finite for every positive integer $m$ such that $m(K_C+B)$ is Cartier.
\end{thm}
\begin{proof}
    We may assume that $C$ is irreducible. After passing to the restricted representations on the subgroup $\Aut(C,B)_F$, which is defined as the subgroup consisting of automorphisms in $\Aut(C,B)$ that fix every irreducible component of $B$, we may assume every element in $\Aut(C,B)$ fixes every component of $B$. To see this, we consider the action of $\Aut(C,B)$ on the finitely many components of $B$, and we know that $\Aut(C,B)_F$ is a normal subgroup of $\Aut(C,B)$ having only finitely many cosets.
    
    If the genus $g(C)>1$, we have that the group $\Aut(C)$ is finite and the theorem follows.  If $g(C)=1$ and $B\neq0$, then $\Aut(C,B)$ is finite. If $g(C)=1$ and $B=0$, we apply Theorem \ref{thm:abelian-variety-pluricanonical-finite}. So we may assume $C\simeq \mathbb{P}^1_t$ with coordinate $t$. Since $(C,B)$ is log canonical with $K_C+B$ nef, it follows that $|\Supp B|>1$. If $|\Supp B|\geq 3$, the group $\Aut(C,B)$ is trivial. If $|\Supp B|=2$, we know that $B$ is reduced since $\deg(K_C+B)\geq0$. We may assume that $B=0+\infty$. In this case, the section $\frac{dt}{t}$ generates $\Oo_C(K_C+B)$. The $m$-th $B$-pluricanonical representation for $(C,B)$ follows.
\end{proof}

\subsection{Cyclic coverings}\label{cyclic covering}
We recall the construction of cyclic coverings in this subsection. The main reference is \cite[2.3]{Kollár_2013}.

Let \( X \) be a normal scheme, \( L \) be a divisorial sheaf on \( X \). We define $L^{[m]}$ to be $(L^{\otimes m})^{**}$. Let \( s \) be a section of \( L^{[m]}\) for some integer \( m > 0 \). Set \( D := (s) \) and \( \Delta := \frac{1}{m}D \). The section \( s \) can be identified with an isomorphism  
\[
\gamma_s: L^{[-m]}(m\Delta) = L^{[-m]}(D) \simeq \mathcal{O}_X.
\]
This in turn defines an algebra structure on  
\[
\mathcal{O}_X + L^{[-1]}([\Delta]) + \cdots + L^{[-(m-1)]}([(m-1)\Delta]),
\]
where, for \( i + j < m \) the multiplication \( L^{[-i]} \times L^{[-j]} \to L^{[-(i+j)]} \) is the tensor product and for \( i + j \geq m \) we compose the tensor product with the isomorphisms  
\[
L^{[-(i+j)]} = L^{[-(i+j-m)]} \otimes L^{[-m]} \xrightarrow{\operatorname{id}\otimes\gamma_s} L^{[-(i+j-m)]} \otimes \mathcal{O}_X = L^{[-(i+j-m)]}.
\]
Then we can define the following cyclic covering over $X$
\[
\tilde{X} = \Spec_X\bigl(\mathcal{O}_X + L^{[-1]}([\Delta]) + \cdots + L^{[-(m-1)]}([(m-1)\Delta])\bigl).
\]

Now we assume $(X,B)$ to be an irreducible pair, with $K_X+B\sim_\Qq0$. Let a positive integer $r$ be the smallest one such that $r(K_X+B)\sim 0$. The integer $r$ is called the \emph{index} of $K_X+B$. We put an additional assumption that $\operatorname{char} k=p\nmid r$.

Fix a non-zero section $\omega \in H^0(X, r(K_X + B))$, then there is an $r$-fold cyclic covering given by $L:=K_X+B$ and $\omega\in H^0(X, r(K_X + B))$
\[
\mu \colon \tilde{X_\omega} = \operatorname{Spec} \left( \bigoplus_{i=0}^{r-1} \mathcal{O}_X(\lfloor i(K_X + B) \rfloor) \right) \to X.
\]
More precisely, the $\mathcal{O}_X$-algebra structure of $\bigoplus_{i=0}^{r-1} \mathcal{O}_X(\lfloor i(K_X + B) \rfloor)$ is given by the natural multiplication
\[
\mathcal{O}_X(\lfloor i(K_X + B) \rfloor) \otimes \mathcal{O}_X(\lfloor j(K_X + B) \rfloor) \to \mathcal{O}_X(\lfloor (i+j)(K_X + B) \rfloor)
\]
if $i+j < r$, and by the multiplication
\begin{align*}
\mathcal{O}_X(\lfloor i(K_X + B) \rfloor) \otimes \mathcal{O}_X(\lfloor j(K_X + B) \rfloor) &\to \mathcal{O}_X(\lfloor (i+j)(K_X + B) \rfloor) \\
&\xrightarrow{\times \omega^{-1}} \mathcal{O}_X(\lfloor (i+j-r)(K_X + B) \rfloor)
\end{align*}
if $i+j \geq r$.

Note that $\tilde{X}_\omega$ is irreducible and reduced by the minimality of $r$ and the assumption that $p\nmid r$, and that $\tilde{X}_\omega$ is normal by the construction (see \cite[2.44.3 and the paragraph above the formula]{Kollár_2013}). Under our assumption $p\nmid r$, the cyclic covering $\mu$ has tame ramification (see \cite[2.41]{Kollár_2013}). There exists a $\mathbb{Q}$-divisor $B_{\tilde{X}_\omega}$ on $\tilde{X}_\omega$ given by
$$K_{\tilde{X}_\omega} + B_{\tilde{X}_\omega} = \mu^*(K_X + B)$$
and $K_{\tilde{X}_\omega} + B_{\tilde{X}_\omega} \sim 0$ by \cite[2.44.6]{Kollár_2013}.


To understand how the $B$-pluricanonical representations of $(X,B)$ interact with the cyclic cover, we give the following lemma.

\begin{lem}[{cf.\cite[Lemma 2.6]{Jiang_2021}}]\label{B bir lifted to cyclic cover}
Let $(X,B)$ be an irreducible projective pair such that $K_X+B\sim_\Qq 0$ with the index of $K_X+B$ being $r$. Assume $\operatorname{char}k=p\nmid r$. Let $\mu:\tilde{X}_{\omega}\to X$ be the $r$-fold cyclic covering given by a non-zero section $\omega\in H^0(X,r(K_X+B))$. Then any automorphism $g\in\Aut(X,B)$ can be lifted to an automorphism $\tilde{g}:(\tilde{X}_{\omega},B_{\tilde{X}_{\omega}})\to(\tilde{X}_{\omega},B_{\tilde{X}_{\omega}})$ satisfying $\mu\circ\tilde{g} = g\circ\mu$, where
$$K_{\tilde{X}_{\omega}}+B_{\tilde{X}_{\omega}}=\mu^*(K_X+B).$$

Furthermore, the finiteness of $A$-pluricanonical representations for $(\tilde{X}_{\omega},B_{\tilde{X}_{\omega}})$ implies the finiteness of $A$-pluricanonical representations for $(X,B)$.
\end{lem}
\begin{proof}
For simplicity, we write $L$ for $K_X+B$. Let
$$g_*(\omega)=c_g\omega$$
for some $c_g\in k^*$ since $H^0(X,rL)$ is one-dimensional. Fix an $r$-th root of $c_g$ and denote it by $c_g^{\frac1r}$. Note that the $r$-fold cyclic covering given by $c_g\omega$ is isomorphic to $\tilde{X}_{\omega}$ by \cite[2.45.2]{Kollár_2013}.

We have the following natural isomorphisms
$$g_*: \Oo_X(\lfloor iL\rfloor)\to \Oo_X(\lfloor iL\rfloor)$$
for any integer $i$. We consider the composition
$$g'_*: \Oo_X(\lfloor iL\rfloor)\xrightarrow{g_*} \Oo_X(\lfloor iL\rfloor)\xrightarrow{\times c_g^{-\frac ir}}\Oo_X(\lfloor iL\rfloor)$$
This induces an isomorphism of the corresponding sheaves of algebras
\[
\begin{tikzcd}[column sep=huge, row sep=large]
\mathcal{O}_{X}\!\left(\left\lfloor i(L)\right\rfloor\right)
\otimes
\mathcal{O}_{X}\!\left(\left\lfloor j(L)\right\rfloor\right)
  \arrow[r]
  \arrow[d,"g'_*\otimes g'_*"']
&
\mathcal{O}_{X}\!\left(\left\lfloor(i+j)(L)\right\rfloor\right)
  \arrow[d,"g'_*"]
\\
\mathcal{O}_{X}\!\left(\left\lfloor i(L)\right\rfloor\right)
\otimes
\mathcal{O}_{X}\!\left(\left\lfloor j(L)\right\rfloor\right)
  \arrow[r]
&
\mathcal{O}_{X}\!\left(\left\lfloor(i+j)(L)\right\rfloor\right)
\end{tikzcd}
\]
for positive integers $i$, $j$ with $i+j<r$;

\[
\begin{tikzcd}[column sep=large, row sep=large]
\Oo_X(\lfloor iL\rfloor) \otimes\Oo_X(\lfloor jL\rfloor)
  \arrow[r]
  \arrow[d,"g'_*\otimes g'_*"']
&
\Oo_X(\lfloor (i+j)L\rfloor)
  \arrow[r,"\times\omega^{-1}"]
  \arrow[d,"g'_*"]
&
\Oo_X(\lfloor (i+j-r)L\rfloor)
  \arrow[d,"g'_*"]
\\
\Oo_X(\lfloor iL\rfloor) \otimes\Oo_X(\lfloor jL\rfloor)
  \arrow[r]
&
\Oo_X(\lfloor (i+j)L\rfloor)
  \arrow[r,"\times \omega^{-1}"']
&
\Oo_X(\lfloor (i+j-r)L\rfloor)
\end{tikzcd}
\]
for positive integers $i$, $j$ with $i+j\geq r$. So this isomorphism gives a birational morphism between coverings
$$\tilde{g}:\tilde{X}_{\omega}\to \tilde{X}_{\omega}.$$
We may consider $g^{-1}$ and then obtain $\tilde{g}^{-1}$ similarly to see the morphism $\tilde{g}$ is in fact an automorphism. To see $\tilde{g}\in\Aut(\tilde{X}_{\omega},B_{\tilde{X}_{\omega}})$, we note that $B_{\tilde{X}_\omega}=\mu^*B-R_\mu$, where $R_\mu$ is the ramification divisor. Given $g^*B=B$ and $\mu\circ\tilde{g} = g\circ\mu$, we conclude that $\tilde{g}$ preserves the divisor $B_{\tilde{X}_\omega}$.

Also, since $K_X+B\sim_\Qq0$, we have
$$H^0(X,r(K_X+B))\simeq H^0(\tilde{X}_{\omega},r(K_{\tilde{X}_{\omega}}+B_{\tilde{X}_{\omega}})).$$
Assume that we have the finiteness of the $A$-pluricanonical representations for $(\tilde{X}_{\omega},B_{\tilde{X}_{\omega}})$, given that
$$\tilde{g}_*(\mu^*\omega)=\mu^*(g_*\omega)=c_g\mu^*(\omega),$$
we know that there exists a positive integer $M$ such that $c_g^M=1$ for any ${g}\in\Aut({X},B_{X})$. Thus, the $A$-pluricanonical representations of $(X,B)$ are also finite.
\end{proof}

\section{Finiteness of $B$-pluricanonical representations}
In this section, we prove the finiteness of $B$-pluricanonical representations. The main theorem of this section is Theorem \ref{main thm: finiteness of pluri for peff lc surf & Cartier ind}.

\subsection{The case of $\kappa(K_S+B)=1$}
In this subsection, we consider log surfaces $(S,B)$ with $\kappa(K_S+B)=1$. Our main theorem in this subsection is Theorem \ref{main thm of kappa one}.

We first investigate the canonical bundle formula in a special case, which will be used in the proof of Theorem \ref{main thm of kappa one} in the case of $\mathbb{P}^1$-bundle. To start with, we recall some definitions.

\begin{definition}[{cf. \cite[Section~7]{PS09}, \cite[Section~6]{cascini2015base}}]\label{div mod}
Let $f\colon X \to Z$ be a surjective proper morphism between two irreducible normal varieties, and suppose that
\[
K_X+D \sim_{\mathbb{Q}} f^*L,
\]
where $L$ is a $\mathbb{Q}$-Cartier $\mathbb{Q}$-divisor on $Z$ and $D$ is a $\Qq$-divisor on $X$ such that $(X,D)$ is log canonical near the generic fiber of $f$, i.e., $(f^{-1}(U),D|_{f^{-1}(U)})$ is log canonical for some Zariski dense open subset $U \subseteq Z$.

We define two divisors $D_{\mathrm{div}}$ and $D_{\mathrm{mod}}$ on $Z$ as follows:
\[
D_{\mathrm{div}}:=\sum_Q(1-c_Q)Q,\qquad D_{\mathrm{mod}}:=L-K_Z-D_{\mathrm{div}},
\]
where the sum runs over all prime Weil divisors $Q \subseteq Z$, and
\[
c_Q=\sup\left\{
c\in\mathbb{R}
\;\middle|\;
(X,D+cf^*Q)
\text{ is log canonical over the generic point }\eta_Q\text{ of }Q
\right\},
\]
so that $K_X+D
\sim_{\mathbb{Q}}
f^*\bigl(K_Z+D_{\mathrm{div}}+D_{\mathrm{mod}}\bigr)$.
\end{definition}

\begin{lemma}\label{canonical bundle formula with trivial moduli}
    Let $f:S\to C$ be a projective surjective morphism, where $S$ is an irreducible normal surface and $C$ is a smooth curve over an algebraically closed field $k$ of characteristic $p >0$. Assume that $f$ is smooth over a Zariski open subset of $C$ with fibers isomorphic to $\mathbb{P}^1$. Let $B$ be a $\Qq$-divisor on $S$ satisfying the following conditions:
    \begin{enumerate}
        \item $(S, B)$ is Kawamata log terminal and $K_S+B\sim_{\Qq,f} 0$;
        \item $\operatorname{char} k=p>\frac{2}{\delta}$, where $\delta$ is the minimal non-zero coefficient of $B^h$, the horizontal part of $B$;
    \end{enumerate}
    Then for a general fiber $F$ of $f$, the pair $(F, B_F)$ defined by $K_F+B_F=(K_S+B)|_F$ is Kawamata log terminal.
    
    Further, suppose that there are infinitely many fibers $(F, \Supp B_F)$ of $f$ that are isomorphic to each other as $n$-pointed stable curves of genus $0$. Then we have $B_{\operatorname{mod}}\sim_\Qq 0$.
\end{lemma}

\begin{proof}
    The lemma follows from the canonical bundle formula in a special case proved in \cite[6.2]{cascini2015base} (see also \cite[Theorem 4.8]{DH16}, \cite[Theorem 8.1]{PS09} and \cite[Lemma 3.1]{Witaszek2021}). By the construction in the cited literature, there is a diagram
\begin{center}
    \begin{tikzcd}
S \arrow[d, "f"'] & \tilde{S} \arrow[d, "g"'] \arrow[l, "\psi"'] \arrow[r, dashed] & \hat{S} \arrow[ld, "\hat f"] \arrow[r, "\hat\phi"] & {\mathcal{U}_{0,n}^*} \arrow[ld, "{g_{0,n}}"] \\
C                 & \tilde{C} \arrow[l, "\psi_0"] \arrow[r, "\phi_0"']             & {\overline{\mathcal{M}}_{0,n}}                     &                                              
\end{tikzcd}
\end{center}
where $\overline{\mathcal{M}}_{0,n}$ is the moduli space of $n$-pointed stable curves of genus $0$, $\mathcal{U}_{0,n}^*$ is a smooth projective variety and the morphism $g_{0,n}:\mathcal{U}_{0,n}^* \to \overline{\mathcal{M}}_{0,n}$ is a $\mathbb{P}^1$-bundle obtained from \cite[6.2]{cascini2015base}, such that
\begin{enumerate}
    \item $\psi_0$ is a tame finite morphism and $\tilde{S}$ is the normalization of the component of $S\times_C\tilde{C}$ dominating $\tilde{C}$. Define $B_{\tilde{S}}$ by $K_{\tilde{S}}+B_{\tilde{S}}=\psi^*(K_S+B)$, then the horizontal components of $B_{\tilde{S}}$ are sections of $g$.
    \item The universal family $f_{0,n}: \overline{\mathcal{U}}_{0,n}\to \overline{\mathcal{M}}_{0,n}$ factors as $f_{0,n}: \overline{\mathcal{U}}_{0,n}\xrightarrow{\sigma}\mathcal{U}_{0,n}^* \xrightarrow{g_{0,n}} \overline{\mathcal{M}}_{0,n}$, where $\sigma$ is a sequence of blowups with smooth centers.
    \item  The morphism $\phi_0$ is induced by the $\mathbb{P}^1$-bundle $g$ and $\hat{S}=\mathcal{U}_{0,n}^*\times_{\overline{\mathcal{M}}_{0,n}}\tilde{C}$.
    \item Denote the sections of $f_{0,n}$, which correspond to the marked points, as $\mathcal{P}_1, \mathcal{P}_2, \cdots, \mathcal{P}_n$, and write the horizontal part of $B_{\tilde{S}}$ as $B^h_{\tilde{S}}=\sum d_i \tilde{B}_i$ with $d_i\neq 0$ for all $i$. Then $(\mathcal{U}_{0,n}^*, \mathcal{D}^*)$ is Kawamata log terminal, where $\mathcal{D}^*:=\sum d_i\sigma_*\mathcal{P}_i$.
    \item    Take the divisors ${B}_{\operatorname{div}}$, ${B}_{\operatorname{mod}}$ on $C$, and the divisors $\tilde{B}_{\operatorname{div}}$, $\tilde{B}_{\operatorname{mod}}$ on $\tilde{C}$ as defined in Definition \ref{div mod}. Then $B_{\operatorname{div}}$ is effective and the following holds
$$K_{\tilde{C}} + \tilde{B}_{\operatorname{div}} = \psi_0^* (K_C+{B}_{\operatorname{div}}),\qquad \tilde{B}_{\operatorname{mod}}\sim_\Qq \psi_0^*{B}_{\operatorname{mod}}.$$
Moreover, there exists a semi-ample $\Qq$-divisor $L$ on $\overline{\mathcal{M}_{0,n}}$ such that
$$\tilde{B}_{\operatorname{mod}}\sim_\Qq\phi_0^*L.$$
\end{enumerate}
The pair $(F, B_F)$ is Kawamata log terminal since it is isomorphic to a fiber of $g_{0,n}$ for a general fiber $F$ of $f$. Now suppose that there are infinitely many fibers of $f$ that are isomorphic as $n$-pointed stable curves of genus $0$. We know that the same holds for the fibration $g$ since $g$ is obtained by fiber product. Therefore the morphism $\phi_0$ is constant, from which we obtain $B_{\operatorname{mod}}\sim_\Qq 0$ and thus $K_S+B\sim_{\Qq}f^*(K_C+B_{\operatorname{div}})$.
\end{proof}

\begin{thm}\label{main thm of kappa one}
    Let $(S,B)$ be a projective irreducible Kawamata log terminal pair of dimension $2$ over an algebraically closed field $k$ of characteristic $p >3$ with $\kappa(K_S+B)=1$. Assume that $K_S+B$ is nef and that $B$ has standard coefficients. Then there exists an integer $m_0$, such that the image of the $A$-pluricanonical representation
    $$\rho_{mm_0}:\Aut(S,B)\rightarrow \Aut\big( H^0(S,mm_0(K_S+B))\big)$$
is finite for every positive integer $m$.
\end{thm}

\begin{proof}
    By abundance for surfaces (see \cite[Theorem 1.2]{abundanceforlcsurf}), the divisor $K_S+B$ is semi-ample with $\kappa(K_S+B)=1$. It follows that we can choose a positive integer $m_0$ and find a fibration between normal varieties $f \colon S \to C$, such that $m_0(K_S+B)=f^*(H)$ for an ample divisor $H$ on $C$, where $C$ is one-dimensional and smooth. In fact, the fibration is obtained as the contraction in the Stein factorization of the morphism
    $$S\to \proj\oplus_{m\geq0}H^0(S, mm_0(K_S+B)).$$
We write $W$ for $\proj \oplus_{m\geq0}H^0(S, mm_0(K_S+B))$. Then we have
$$S\xrightarrow{f} C\xrightarrow{r} W$$
and $H=r^*(\Oo_W(1))$.

    Since elements in $\Aut(S,B)$ preserve sections in $H^0(S, mm_0(K_S+B))$, we know that the group $\operatorname{Aut}(S,B)$ acts on $S$ preserving this fibration $f$. More precisely, every $g\in \operatorname{Aut}(S,B)$ induces an automorphism of $H^0(S, mm_0(K_S+B))$ for every integer $m$, and thus an automorphism of the algebra
$$R:=\oplus_{m\geq0} H^0(S, mm_0(K_S+B)),$$
which we denote by $\rho(g)$. By the construction of $C$, the induced automorphism of $W$ induces an automorphism of $C$, which we denote by $l(g)$.

To conclude, we have the following commutative diagram:
\[
\begin{tikzcd}
(S, B) \arrow[r, "g"] \arrow[d, "f"] & (S, B) \arrow[d, "f"] \\
C \arrow[r, "l(g)"] & C.
\end{tikzcd}
\]
Also, we can construct the following group homomorphism
$$\Aut(S,B)\xrightarrow{\rho}\Aut R\xrightarrow{i}\Aut C$$
with the composition $i\circ\rho=l$. And to prove the finiteness of the image $\rho_{mm_0}(\Aut(S,B))$, it suffices to prove the finiteness of the image $\rho(\Aut(S,B))$. For that purpose, we claim the following,

\begin{itemize}
    \item Claim \ref{claim:kernel}: The kernel of the map $i|_{\operatorname{Im}\rho}$, which is the map $i$ restricted on the subgroup $\operatorname{Im}\rho$ of $\Aut R$, is finite ;
    \item Claim \ref{claim:image}: The image of $l=i\circ\rho$ is finite.
\end{itemize}

\begin{claim}\label{claim:kernel}
The kernel of the restriction $i|_{\operatorname{Im}\rho}$ is finite;
\end{claim}

\begin{proof}[Proof of Claim \ref{claim:kernel}]
For an element $g\in\operatorname{ker}l$, the automorphism $\rho(g)$ acts on the linear spaces $H^0(S,m_0(K_S+B))$ by scalar multiplication since $\rho(g)$ acts on $\proj R$ as identity. Thus, to prove the claim, we need to show that there exists a positive integer $M$ such that for any element $g\in\operatorname{ker}l$, and any non-zero element $\alpha\in H^0(S,mm_0(K_S+B))$ such that $g^*\alpha=a\alpha$, we have $a^M=1$.

Take the generic fiber $F$, which is one-dimensional, normal and Gorenstein by the fact that $S$ is normal. Also, $F$ is of Calabi-Yau type since $f$ is $K_S+B$-trivial. Suppose $F$ is not geometrically normal. Consider the normalization of the reduced subscheme $F_0$ of the scheme $\bar F:=F\otimes_{k(C)}\Spec\ \overline{k(C)}$, which we denote by $\pi: T\to F_0$. We also denote the composition morphism $T\xrightarrow{\pi} F_0\hookrightarrow \bar F\to F$ as $\phi$. Then by \cite[Theorem 1.1]{PW17}, we have
$$K_T+(p-1)D\sim \phi^*K_{F},$$
where $D$ is an integral effective divisor (not $0$ if and only if $F$ is not geometrically normal). Thus, we have
$$-2\leq\deg K_T\leq -\deg (p-1)D\leq -p+1, $$
which is impossible given our assumption $\operatorname{char} k=p>3$.

So $F$ is geometrically normal and hence smooth. By adjunction on $F$, we obtain a pair $(F,B_F)$ with $K_F+B_F\sim_\Qq 0$. Note that $l(g)$ acts on $C$ as the identity. Hence, the morphism $g$ preserves all fibers of $f$ and induces an automorphism of the curve $(F,B_F)$. Therefore, if an element $\alpha\in H^0(S,mm_0(K_S+B))$ is non-zero and $g^*\alpha=a\alpha$, then after restricting to $F$, we get
\[
(g|_F)^*(\alpha|_F)=a(\alpha|_F),
\]
where $\alpha|_F\neq0$ since $F$ is the generic fiber. Thus, to bound the order of $a$, it is enough to prove the finiteness of $A$-pluricanonical representations for the curve $(F,B_F)$. We conclude by Theorem \ref{1-dim subklt fin rep}.
\end{proof}

\begin{claim}\label{claim:image}
The image of $l$ is finite.
\end{claim}

\begin{proof}[Proof of Claim \ref{claim:image}]
Take the minimal resolution $g:S'\to S$. We have
$$B':=g^*(K_S+B)-K_{S'}\geq0.$$
Since $g$ is a minimal resolution, every automorphism of $S$ induces an automorphism of the smooth surface $S'$. Similar to Remark \ref{Main Thm equi under B-bir}, by replacing $(S,B)$ with $(S',B')$, we may assume that $S$ is smooth. Passing the subgroup of $\Aut(S,B)$ consisting of elements fixing all irreducible components of reducible fibers of $f$, which is of finite index, we may assume $\Aut(S,B)$ fixes all $(-1)$-curves in the fibers of $f$. By contracting $(-1)$-curves in the fibers of $f$ (which is $\Aut(S,B)$-invariant), we may assume $f: S\to C$ is a minimal elliptic fibration or a $\mathbb{P}^1$-bundle. Since \(\operatorname{char}k>3\), the quasi-elliptic case does not occur; indeed, quasi-elliptic fibrations exist only in characteristics \(2\) and \(3\), by \cite[Corollary~4.2]{liedtke2013algebraic}. Assume the image of $l$ is infinite for contradiction. Thus $C$ has infinite automorphisms that preserve the polarization given by $H$. Thus $C\simeq \mathbb{P}^1$.

First, we suppose that $f$ is a $\mathbb{P}^1$-bundle. From the proof of the Claim \ref{claim:kernel}, we know the generic fiber of $f$ is geometrically smooth. Therefore $f$ is smooth over a Zariski open subset of $C$. Also, since $B$ has standard coefficients and $p> 3$ by our assumption, we have that $\operatorname{char}k=p>\frac{2}{\delta}$, where $\delta$ is the minimum non-zero coefficient of $B^h$. With our assumption that $\operatorname{Im}(l)$ is infinite, we know that there are infinitely many fibers $(F, \Supp B_F)$ of $f$ that are isomorphic to each other as $n$-pointed stable curves of genus zero. By Lemma \ref{canonical bundle formula with trivial moduli}, we obtain $K_S+B\sim_\Qq f^*(K_C+B_{\operatorname{div}})$ with $(C,B_{\operatorname{div}})$ log canonical and $K_C+B_{\operatorname{div}}$ ample. By the construction of $B_{\operatorname{div}}$, the map $l:\Aut(S,B)\to \Aut C$ factors through $\Aut(C, B_{\operatorname{div}})\to \Aut C$. Since $(C,B_{\operatorname{div}})$ is log canonical with $K_C+B_{\operatorname{div}}$ ample, the group $\Aut(C, B_{\operatorname{div}})$ is finite (see the proof of Theorem \ref{1-dim subklt fin rep}), contradicting our assumption that $\operatorname{Im}(l)$ is infinite.

If the fibration $f:S\to C$ is a minimal elliptic fibration, $B^h=0$. By \cite[Theorem 7.15 (and 7.3.4)]{Gu2020SurfaceFW}, we only come to several cases below. We first discuss the possible locus of the boundary divisor. Note that the surface $S$ in each case below has a nontrivial action by either the linear algebraic group $\ga$ or $\gm$ lifted from the base. Therefore, by \cite[Corollary 3]{Brion2021}, $\mathbf{Aut}(S)$ is an algebraic group, and thus so is the closed subgroup subscheme $\mathbf{Aut}(S,B)_F:=
\bigcap_i\operatorname{Stab}_{\mathbf{Aut}(S)}(B_i)
$, the stabilizer subgroup of $B$ (see \cite[II.1.3.6]{DemazureGabriel1970} or \cite[VIB.6.2.4]{SGA3}). Here $B_i$ are components of $B$ with non-zero coefficients. Therefore the image of the morphism between group schemes
$$\mathbf{Aut}(S,B)_F\to \mathbf{Aut}(C)$$
is closed by \cite[Proposition 1.68]{milne2017algebraicgroups}, and has infinite closed points by assumption. We denote the neutral component of the image by $G$, which is a nontrivial linear algebraic group ($G$ preserves the polarization on $C$ given by $H$). By \cite[Lemma 7.2]{Chen_2024}, $G$ contains either $\ga$ or $\gm$, such that $G$ acts on $C\simeq\mathbb{P}^1$ additively or multiplicatively respectively. Since $B=B^v$ is invariant under the action of $\Aut(S,B)_F$, it follows that it can be supported only over the point ${0}$ and $\infty$. Here, $B^v$ is defined to be the vertical part of $B$ with respect to $f$. We will deduce that $K_S+B\leq 0$ in each case, contradicting the fact that $K_S+B=K_S+B^v$ is a pull back of an ample $\Qq$-divisor on $C$ along $f$.

We define Weil divisors $F_0:=f^{*}(0)$, $F_\infty:=f^{*}(\infty)$, and Weil divisors $L_0$, $L_\infty$ to be the reduced divisors having the same support as $F_0$, $F_\infty$ respectively.

\begin{enumerate}[label=(\arabic*), labelsep=0.5em]
    \item The fibration $f$ is the projection
    $$f=p_2:S:=E\times\mathbb{P}^1\to \mathbb{P}^1$$
        for an elliptic curve $E$.
        
        In this case, $K_S=f^*(K_{\mathbb{P}^1})$. Since $(S,B)$ is log canonical, we have $K_S+B\sim_\Qq f^*(-2+\deg f_*(B))\leq 0$.

    \item There is an elliptic curve \( E \) over \( k \), a non-trivial subgroup scheme \( \Lambda \subseteq E \) along with an embedding \( \Lambda \subseteq \mathbb{G}_a \), so \( \Lambda \) acts on \( \mathbb{P}^1 \) additively. The fibration \( f \) is as below:
\[
f = p_2 : S := (E \times \mathbb{P}^1)/\Lambda \to \mathbb{P}^1/\Lambda.
\]
Here \( \Lambda \) acts on \( E \times_k \mathbb{P}^1_k \) diagonally. The complete classification of \( (E, \Lambda) \) is as follows:
\begin{itemize}
    \item for ordinary \( E \) over \( k \), we take \( \Lambda = E[p](k) \);
    \item for supersingular \( E \), we take \( \operatorname{Ker}(F_{E/k} : E \to E^{(p)}) \simeq \alpha_p \).
\end{itemize}

By \cite[Example 4.7, 4.8]{katsura1984elliptic}, we have that $K_S=f^*(-1)+(p-2)L_\infty$, where $L_\infty$ is the only singular fiber. And note that the $\gm$ action on $\mathbb{P}^1/\Lambda$ does not lift to $S$. It follows that $B$ can be supported only over $\infty$. Since $(S,B)$ is log canonical, $L_\infty-B\geq 0$, and thus
$$K_S+B<K_S+B - (f^*(-1)+(p-2)L_\infty+L_\infty)\leq 0.$$

    \item There is an elliptic curve \( E \), a subgroup scheme \( \Lambda \subseteq E \) along with an embedding \( \Lambda \subseteq \mathbb{G}_m \), so \( \Lambda \) acts on \( \mathbb{P}^1_k \) multiplicatively. The fibration \( f \) is as below:
\[
f = p_2 : S := (E \times \mathbb{P}^1) / \Lambda \to C = \mathbb{P}^1 / \Lambda.
\]
Here \( \Lambda \) acts on \( E \times \mathbb{P}^1 \) diagonally. The complete classification of \( (E, \Lambda) \) is as follows:
\begin{itemize}
    \item for any \( E \), we can take \( \Lambda \subseteq E(k) \) to be a non-trivial cyclic subgroup of order prime to char \( k \);
    \item when \( E \) is not supersingular, we can take \( \Lambda \) to be either \( \operatorname{Ker}(F^{(r)} : E \to E^{(p^r)}) \simeq \mu_{p^r} \) or its product with a non-trivial cyclic subgroup of \( E(k) \) of order prime to \(\operatorname{char} k \).
\end{itemize}

Set $n:=\dim_k\mathcal O_\Lambda$. Then $F_0=nL_0, F_\infty=nL_\infty$. Denote the coordinate of $\mathbb{P}^1$ to be $t$. Consider the projection $S\to C'=E/\Lambda$, with two sections $L_0, L_\infty$, we have that $\omega_{S/C'}(L_0+L_\infty)\simeq\mathcal O_S$ and thus
$$K_S+L_0+L_\infty\simeq \Oo_S.$$
Since $(S,B)$ is log canonical, $(\frac{1}{n}F_0+\frac{1}{n}F_\infty)-B\geq 0$, and thus
$$K_S+B=K_S+B - \big((-\frac{1}{n}F_0-\frac{1}{n}F_\infty)+L_0+L_\infty\big)\leq 0.$$

\item There is an elliptic curve \( E \), a non-trivial cyclic group \( \Gamma \subseteq \operatorname{Aut}_k^{\text{gp}}(E) \) of order \( n \) prime to \( \operatorname{char}k \). Fix an arbitrary isomorphism \( \Gamma \simeq \mu_n \), then \( \Gamma \) also acts on \( \mathbb{P}^1 \) multiplicatively. The fibration \( f \) is as below:
\begin{center}
\begin{tikzcd}[column sep=10em, row sep=4em]
S
  \arrow[r,dashed, "\textit{relatively minimal regular model}",
    "\vartheta"']
  \arrow[d, "f"']
&
(E\times\mathbb P^1)/\Gamma
  \arrow[d, "p_2"]
\\
C
  \arrow[r, equal]
&
\mathbb P^1/\Gamma .
\end{tikzcd}
\end{center}
Here \( \Gamma \) acts on \( E \times \mathbb{P}^1 \) diagonally. The complete classification of \( (E, \Gamma) \) is as follows:
\begin{itemize}
    \item \( p \neq 2,\; \Gamma = \{\pm 1\} \). The singular fibres at \( 0, \infty \) are of type \( I_0^* \);
    \item \( p \neq 3,\; E : y^2 - y = x^3 - 1 \) (\( j(E) = 0 \)), \( \Gamma = \langle \gamma \rangle \) with \( \gamma : (x, y) \mapsto (\zeta_3 x, y) \). The singular fibres at \( 0, \infty \) are of type \( IV \) and \( IV^* \) (or vice versa) respectively;
    \item \( p \neq 2,\; E : y^2 = x^3 - x \) (\( j(E) = 1728 \)), \( \Gamma = \langle \gamma \rangle \) with \( \gamma : (x, y) \mapsto (-x, \zeta_4 y) \). The singular fibres at \( 0, \infty \) are of type \( III \) and \( III^* \) (or vice versa) respectively.
    \item $p\neq2,3$, $E:y^2=x^3-1$ ($j(E)=0$), \( \Gamma = \langle \gamma \rangle \) with \( \gamma : (x, y) \mapsto (\zeta_3x, - y) \). The singular fibres at \( 0, \infty \) are of type \( II \) and \( II^* \) (or vice versa) respectively.
\end{itemize}
    In particular, the fibration has no multiple fibers. Hence \cite[(1.7), (1.8)]{katsura1984elliptic} implies that $ \operatorname{Tor}(R^1f_*\mathcal O_S)=0$. Moreover, $S$ is rational. Indeed, projection to $E$ induces a rational fibration $S\dashrightarrow E/\Gamma\simeq\mathbb P^1$ whose generic fiber is $\mathbb P^1$. Therefore we conclude $\chi(\mathcal O_S)=1$. By \cite[(1.3), (1.4)]{katsura1984elliptic}, we have the canonical bundle formula
\[
\begin{aligned}
    \omega_S
    &\simeq
    f^*\mathcal O_{\mathbb P^1}(-1),
\end{aligned}
\]
implying $B$ is relatively $\Qq$-linearly trivial. The log canonical threshold of these types of fibers are listed in the following table (see \cite[\textbf{Degenerate fibers of (quasi-)elliptic fibrations} in Section 4]{liedtke2013algebraic}, \cite[Table 1]{bejleri2026moduli})
\[
\begin{array}{c|ccccccc}
 &I_0^*&II&II^*&III&III^*&IV&IV^*\\ \hline
\operatorname{lct}(S;F)
&\frac12&\frac56&\frac16&\frac34&\frac14&
\frac23&\frac13
\end{array}
\]
Therefore, if $B=b_0F_0+b_\infty F_\infty$, we have $b_0+b_\infty\leq \lct(S,F_0)+\lct (S,F_\infty)=1$. Thus, $K_S+B\leq 0$.
\end{enumerate}

To summarize, we exclude all these cases. Otherwise, $m_0(K_S+B)$ can not be semi-ample with $\kappa(K_S+B)=1$, which is a contradiction. Therefore, the image $\operatorname{Im}(l)$ is finite.
\end{proof}
\noindent \textit{Proof of Theorem \ref{main thm of kappa one} continued.} Given the two claims, we can prove the finiteness of the image $\rho(\operatorname{Aut}(S,B))$, and thus the finiteness of the image $\rho_{mm_0}$.
\end{proof}
\begin{rmk}
    In the proof in the case of $\mathbb{P}^1$-bundle above, given $B_{\operatorname{mod}}\sim_\Qq 0$, another way to deduce contradiction is as follows:
    
    Note that $C\simeq \mathbb{P}^1$ under our assumption that $\operatorname{Im}(l)$ is infinite, and that canonical bundle formula gives $K_S+B\sim_{\Qq} f^*(K_C+B_{\operatorname{div}})$ with $(C,B_{\operatorname{div}})$ log canonical. Similar to the proof in the case of elliptic fibration, we conclude that $\Supp B_{\operatorname{div}}\subset\{0,\infty\}$ (Note that $B$ is big in this case and thus $\mathbf{Aut}(S,B)$ is a linear algebraic group). Thus we have that $K_S+B\le f^*(K_C+\{0\}+\{\infty\})=0$, which is a contradiction.
\end{rmk}

\subsection{The case of $\kappa(K_S+B)=0$}
In this subsection, we study the finiteness of $A$-pluricanonical representations for log Calabi-Yau surfaces using cyclic covering. The main result of this section is Theorem \ref{main thm of cy surf}.

First, we prove the finiteness of $A$-pluricanonical representations for abelian varieties in characteristic $p>0$.

\begin{thm}\label{thm:abelian-variety-pluricanonical-finite}
Let $A$ be an irreducible abelian variety over an algebraically closed field $k$ of characteristic $p>0$. Then for every integer $m\geq 1$, the image of the $A$-pluricanonical representation
\[
\rho_m:\Aut(A)\longrightarrow \Aut\big( H^0(A,mK_A)\big)
\]
is finite.
\end{thm}

\begin{proof}
Put $V=H^0(A,\Omega_A^1)$. Global one-forms on an abelian variety are translation-invariant (see \cite[Chapter I, Corollary 1.2; Chapter IV, Proposition 6.4]{Milne1986}). Thus, $V$ is isomorphic to the cotangent space at the origin of $A$.

Consider the natural map given by pullback
\[
\Delta:\operatorname{End}(A)\longrightarrow \operatorname{End}_k(V),\qquad
\alpha\mapsto \alpha^*|_V .
\]
We claim that this map is additive (although it reverses composition).

Let $\alpha,\beta\in \operatorname{End}(A)$, let $m_A:A\times A\to A$ be the group law, and let $\delta_A:A\to A\times A$ be the diagonal morphism. Then the sum endomorphism $\alpha+\beta$ is the composition
\[
A\xrightarrow{\delta_A}A\times A
\xrightarrow{\alpha\times\beta}A\times A
\xrightarrow{m_A}A.
\]
For every element $\eta\in V$, translation invariance gives
\[
m_A^*\eta=\operatorname{pr}_1^*\eta+\operatorname{pr}_2^*\eta
\]
in $H^0(A\times A,\Omega^1_{A\times A})$. In fact, denote the translation induced by a point $x$ on $A$ by $t_x$, then for any point $(x,y)$ on $A\times A$ and any tangent vector $(v,w)$ at this point with the identification $T_{(x,y)}(A\times A)=T_xA\oplus T_yA$, we have
$$d(m_A)_{(x,y)}(v,w)=d(t_y)_x(v)+d(t_x)_y(w),$$
and thus
\[
\begin{aligned}
(m_A^*\eta)_{(x,y)}(v,w)
&=\eta_{x+y}\bigl(d(m_A)_{(x,y)}(v,w)\bigr)\\
&=\eta_{x+y}\bigl(d(t_y)_x(v)\bigr)
  +\eta_{x+y}\bigl(d(t_x)_y(w)\bigr)\\
&=\eta_x(v)+\eta_y(w)\\
&=(\operatorname{pr}_1^*\eta+\operatorname{pr}_2^*\eta)_{(x,y)}(v,w).
\end{aligned}
\]
Therefore
\[
\begin{aligned}
(\alpha+\beta)^*\eta
&=\delta_A^*(\alpha\times\beta)^*m_A^*\eta \\
&=\delta_A^*(\alpha\times\beta)^*(\operatorname{pr}_1^*\eta+\operatorname{pr}_2^*\eta) \\
&=\delta_A^*(\operatorname{pr}_1^*\alpha^*\eta+\operatorname{pr}_2^*\beta^*\eta) \\
&=\alpha^*\eta+\beta^*\eta,
\end{aligned}
\]
where the last equality follows from $\operatorname{pr}_1\circ\delta_A=\operatorname{id}_A$ and $\operatorname{pr}_2\circ\delta_A=\operatorname{id}_A$. Thus $\Delta$ is additive.

By \cite[Chapter I, Theorem 10.15]{Milne1986}, $\operatorname{End}(A)$ is a free abelian group of finite rank. Choose a $\mathbb Z$-basis $e_1,\ldots,e_r$ of $\operatorname{End}(A)$. If $\alpha=\sum_i n_i e_i$, then
\[
\Delta(\alpha)=\sum_i n_i\Delta(e_i).
\]
Over the field $k$, this expression depends only on the classes of the integers $n_i$ modulo $p$. Hence $\Delta(\operatorname{End}(A))$ is finite, and so is $\Delta(\Aut_{\mathrm{gp}}(A))$.

Finally, since every automorphism of $A$ as a variety is a translation composed with a group automorphism (see \cite[Chapter I, Corollary 1.2]{Milne1986}), and translations act trivially on global one-forms (see \cite[Chapter IV, Proposition 6.4]{Milne1986}), it follows that the image of the following pullback map
$$\Aut(A)\to\Aut(V),\qquad
\alpha\mapsto \alpha^*|_V $$
is finite. Given
\[
H^0(A,\omega_A^{\otimes m})\simeq \left(\bigwedge^{\dim A}V\right)^{\otimes m},
\]
the induced action of $\Aut(A)$ on $H^0(A,\omega_A^{\otimes m})$ has finite image. This proves the theorem.
\end{proof}

The following theorem gives the finiteness of $A$-pluricanonical representations for smooth $K3$ surfaces.

\begin{thm}[{cf. \cite[Lemma 2.1]{Matsumoto_2023}, \cite[Proposition 3.5]{jang2016representationsautomorphismgroupsfrobenius}, \cite[Theorem 2.1]{Nygaard1980HigherDR}}]\label{main thm for smooth K3 surface}
    Let $X$ be a $K3$ surface over a field $k$. Then $\rho_m(\Aut(X))$ is a finite cyclic group.
\end{thm}

Now we prove the finiteness of $A$-pluricanonical representations for log Calabi-Yau surfaces under some further assumptions as follows.

\begin{thm}\label{main thm of cy surf}
    Let $(S,B)$ be a projective Kawamata log terminal pair of dimension $2$ over an algebraically closed field $k$ of characteristic $p >42$ such that $K_S+B\sim_{\mathbb Q}0$. Assume that $B$ has standard coefficients. Then there exists a positive integer $m_0$ such that the $A$-pluricanonical representation
$$\rho_{mm_0}:\Aut(S,B)\rightarrow \Aut\big( H^0(S,mm_0(K_S+B))\big)$$
is finite for every positive integer $m$.
\end{thm}
\begin{proof}
We may assume $S$ is irreducible. Let $n$ be the index of $K_S+B$. By the bound of global indices for two-dimensional log Calabi--Yau pairs proved in \cite{BWX26indices}, we know that $p\nmid n$.

Choose a non-zero section $\omega\in H^0(S,n(K_S+B))$,
and let
\[
        \pi:\tilde S_\omega\to S
\]
be the corresponding $n$-fold cyclic index-one cover. By construction, we have
$$K_{\tilde S_\omega}=\pi^*(K_S+B)\sim0,$$
since the coefficients of $B$ are standard. Given that the cover is tame, it follows that the surface $\tilde S_\omega$ is Gorenstein and Kawamata log terminal (see \cite[Corollary 2.43, Definition 2.49, Proposition 2.50]{Kollár_2013}). Thus, $\tilde S_\omega$ has Du Val singularities. Let $\mu:S'\to \tilde S_\omega$ be the minimal resolution of $\tilde{S}_\omega$. We have
\[
        K_{S'}=\mu^*K_{\tilde S}\sim0.
\]
Therefore, $S'$ is a smooth projective minimal surface with linearly trivial canonical divisor. By the classification of smooth projective minimal surfaces of Kodaira dimension zero in characteristic $p>3$ (see \cite[Chapter 10]{buadescu2001algebraic}), the surface $S'$ is either a $\operatorname{K3}$ surface or an abelian surface.

By Lemma \ref{B bir lifted to cyclic cover}, Remark \ref{Main Thm equi under B-bir}, we may reduce the question to $(S',0)$, and conclude by Theorem \ref{thm:abelian-variety-pluricanonical-finite} and
Theorem \ref{main thm for smooth K3 surface}.
\end{proof}

\subsection{Proof of Theorem \ref{main thm: finiteness of pluri for peff lc surf & Cartier ind}}
We give the proof of the finiteness of $B$-pluricanonical representations for log canonical surfaces in this subsection. To be precise, we prove the following result.

\begin{thm}\label{main thm: finiteness of pluri for peff lc surf & Cartier ind}
Let $(S,B)$ be a projective log canonical pair of dimension $2$ over an algebraically closed field $k$ of characteristic $p >42$. Assume that $K_S+B$ is nef. Then the image of the $B$-pluricanonical representation
    \[
        \rho_{m}:\Bir(S,B)\longrightarrow
        \Aut\big( H^0(S,m(K_S+B))\big)
    \]
is finite for every positive integer $m$ such that $m(K_S+B)$ is Cartier.
\end{thm}

First, we consider Kawamata log terminal pairs $(S,B)$ and prove Theorem \ref{main thm: finite pluri repn for Kawamata log terminal}. The following lemma shows that we may assume $K_S+B$ is semi-ample by running a $(K_S+B)$-MMP.

\begin{lem}\label{emmp}
Let $(S,B)$ be a projective irreducible log canonical surface. Assume that $S$ is $\Qq$-factorial and $K_S+B$ is pseudo-effective. Then there exists a $\Aut(S,B)$-equivariant birational contraction
\[
        \pi:S\to S_{\min}
\]
such that, if $B_{\min}:=\pi_*B$, then $(S_{\min},B_{\min})$ is log canonical and $K_{S_{\min}}+B_{\min}$ is nef. The equivariance means that, every automorphism in $\Aut(S,B)$ descends to the minimal model: for every $g\in \Aut(S,B)$ there exists a unique automorphism
\[
        g_{\min}\in \Aut(S_{\min},B_{\min})
\]
such that
\[
        \pi\circ g=g_{\min}\circ \pi .
\]
Furthermore, the finiteness of $A$-pluricanonical representations for $(S_{\min},B_{\min})$ implies the finiteness of $A$-pluricanonical representations for \((S,B)\).
\end{lem}

\begin{proof}
Write the Zariski decomposition
\[
        K_S+B=P+N=P+\sum_{i=1}^{q}a_iC_i,
\]
where $P$ is a nef $\Qq$-divisor, $N\geq 0$, and if $N > 0$, then $(P \cdot C_i) = 0$, $\forall i=1,2,\cdots,q$ and the intersection matrix of $\{C_i\}_{i=1}^{q}$ is negative definite.

By running $(K_S+B)$-MMP, we contract $(K_S+B)$-negative curves. If $N=0$, implying that $K_S+B$ is nef, the assertion is immediate. Otherwise, we have $N\neq 0$. If a curve $C$ is contracted at the first step of MMP, which we denote as $f:S\to S_1$ , then $(K_S+B)\cdot C<0$, implying $N\cdot C<0$ since $P$ is nef. Thus, $C=C_i$ for some $i$. We may assume $i=q$. Then we have
\begin{equation}
K_{S_1}+B_1:=f_*(K_S+B)
=f_*(P+N)
=f_*P+\biggl(\sum_{i=1}^{q-1} a_i f_*C_i\biggr).
\tag{$\star$}
\end{equation}
For $i=1,\cdots, q-1$, we have $f_*C_i\neq 0$ since $C_i$ cannot be contracted (otherwise, the curve $C_i$ is numerically equivalent to a scalar multiple of $C_q$, contradicting the fact that the intersection matrix of $\{C_i\}_{i=1}^{q}$ is negative definite).

By projection formula, the divisor $f_*P$ is nef, and has zero intersection numbers with the curves $f_*C_i$ for $i=1,\cdots, q-1$. And also, the intersection matrix of $\{f_*C_i\}_{i=1}^{q-1}$ is equal to that of
$$\{f^*f_*C_i=C_i-\frac{C_i\cdot C_q}{C_q^2}C_q\}_{i=1}^{q-1}$$
which is negative definite. By uniqueness of Zariski decomposition, we conclude that $({\star})$ gives the Zariski decomposition of $K_{S_1}+B_1$ with the positive part being $f_*P$ and the negative part being $\sum_{i=1}^{q-1} a_if_*C_i$. Similarly, every $(K_{S_1}+B_1)$-negative curve is contained in $\{f_*C_i\}_{i=1}^{q-1}$.

Continue this process, we know that the resulting morphism of the MMP, $\pi:S\to S_{\min}$ contracts exactly all the curves in $\{C_i\}_{i=1}^{q}$.

Since every element of $\Aut(S,B)$ preserves $K_S+B$, the uniqueness of Zariski decomposition implies that $\Aut(S,B)$ preserves the divisor $N$. So every $g\in\Aut(S,B)$ maps any exceptional divisor of $\pi$ to an exceptional divisor of $\pi$. Hence, there is a unique automorphism $g_{\min}$ of $S_{\min}$ such that
\[
        \pi\circ g=g_{\min}\circ\pi .
\]
Furthermore, since $g$ preserves $B$, we have
\[
        (g_{\min})_*B_{\min}
        =(g_{\min})_*\pi_*B
        =\pi_*g_*B
        =\pi_*B
        =B_{\min},
\]
and therefore $g_{\min}\in \Aut(S_{\min},B_{\min})$. This proves the automorphism descent statement and gives the homomorphism
$$h:\Aut(S,B)\to \Aut(S_{\min},B_{\min}).$$
Finally, we have a natural isomorphism between global sections
$$\pi_*: H^0(S,m(K_S+B))\to H^0(S_{\min},m(K_{S_{\min}} + B_{\min}))$$
for any integer $m$ that is sufficiently divisible. Also, for any section $\omega\in H^0(S,m(K_S+B))$, we have
$$\pi_{*}\big(\rho_{m}(g)(\omega)\big)=\pi_{*}g_*\omega={g_{\min}}_*\pi_*\omega.$$
Hence, if we denote $V_m$ to be $H^0(S,m(K_S+B))$ and $H^0(S_{\min},m(K_{S_{\min}}+B_{\min}))$ under the identification given by $\pi_*$, then the $A$-pluricanonical representation $\rho_m$ of $(S,B)$
$$\rho_m:\Aut(S,B)\to \Aut(V_m)$$
factors through
$$h:\Aut(S,B)\to \Aut(S_{\min},B_{\min}).$$
So if the $A$-pluricanonical representations of $(S_{\min},B_{\min})$ is finite, then the $A$-pluricanonical representations of $(S,B)$ is also finite. This proves the lemma.
\end{proof}

\begin{thm}\label{main thm: finite pluri repn for Kawamata log terminal}
    Let $(S,B)$ be a projective Kawamata log terminal pair of dimension $2$ over an algebraically closed field $k$ of characteristic $p >42$. Then there exists a positive integer $m_0$ such that the image of the representation
    $$\rho_{mm_0}:\Bir(S,B)\rightarrow\Aut\big( H^0(S,mm_0(K_S+B))\big)$$
is finite for every positive integer $m$.
\end{thm}
\begin{proof}
We may assume $S$ is irreducible. Take a terminal model $(S',B')$ of $(S,B)$, whose existence is guaranteed by \cite[Lemma 2.22]{haconxu2015}. Equivalently speaking, there exists a birational morphism $f:S'\to S$ and a divisor $B'$ on $S'$ such that the pair $(S',B')$ is terminal and $K_{S'}+B'=f^*(K_S+B)$. For any $B$-birational map $\sigma:S\dashrightarrow S$, the composite $\sigma':=f^{-1}\circ\sigma\circ f$ is a birational self-map of $S'$. Moreover, it is an automorphism. To see this, we consider the following diagram
\[
\xymatrix{
& W \ar[dl]_{g_1} \ar[dr]^{g_2} & \\
S' \ar@{-->}[rr]^{\sigma'} \ar[d]_f & & S' \ar[d]^f \\
S \ar@{-->}[rr]^{\sigma} & & S
}
\]
where $W$ is a common resolution of $f$ and $\sigma'$. Then we have
\[
g_1^*(K_{S'}+B')=(f\circ g_1)^*(K_S+B)=(f\circ g_2)^*(K_S+B)=g_2^*(K_{S'}+B').
\]
Now we write
\[
K_W+B_W=g_1^*(K_{S'}+B')=g_2^*(K_{S'}+B').
\]
Since $(S',B')$ is terminal, it follows that any exceptional divisor of $g_1,g_2$ has positive discrepancy, and any non-exceptional divisor of $g_1,g_2$ has non-positive discrepancy. Thus, the curves on $W$ that are contracted by $g_2$ are exactly those that have negative coefficients in $B_W$, and hence they are also contracted by $g_1$. This shows that $\sigma'$ is an automorphism.

By Remark \ref{Main Thm equi under B-bir}, we may replace $(S,B)$ by $(S',B')$, and assume that $(S,B)$ is terminal and that $\Bir(S,B)=\Aut(S,B)$. Also, we may assume $K_S+B$ is pseudo-effective.

Let $h:\widetilde{S}\to S$ be the minimal log resolution of the pair $(S,B)$. Write
\[
        K_{\widetilde{S}}+h_*^{-1}B=h^*(K_S+B)+\sum_i a_iE_i,
\]
where the $E_i$ are $h$-exceptional divisors. Since $(S,B)$ is terminal, $a_i>0$ for all $i$.

For a positive integer $N$, define
\[
        D_N:=h_*^{-1}B+\sum_i \frac1{N} E_i,
        \qquad
        \Delta_{N}:=\left(1-\frac1{N}\right)\lceil D_N\rceil .
\]

The log smooth pair $(\widetilde{S},\Delta_{N})$ is also Kawamata log terminal, and $\Delta_{N}$ has standard coefficients.

Fix a positive integer $N$ that is sufficiently divisible and large (in particular, $N>2$), so that $N(K_S+B)$ is Cartier and the coefficients of $Nh_*^{-1}B$ are integers. As $(S,B)$ is Kawamata log terminal, the coefficients of $Nh_*^{-1}B$ are at most $N-1$. Hence
\[
        D_N\leq (1-1/N)\lceil D_N\rceil=\Delta_N.
\]
And for any positive integer $m$, we have
\[
        mN(K_{\widetilde{S}}+D_N)=h^*(mN(K_S+B))+\sum_i(mNa_i+m)E_i .
\]
Suppose now we can prove the finiteness of $A$-pluricanonical representations for $(\widetilde{S},\Delta_{N})$. Namely, there exists an integer $m_0'$ such that the image of the representation
$$\Aut(\widetilde{S},\Delta_{N})\rightarrow \Aut(H^0(\widetilde{S},mm_0'(K_{\widetilde{S}}+\Delta_{N})))$$
is finite for every $m$.

Note that by the uniqueness of the minimal log resolution of a surface pair, every automorphism of $(S,B)$ lifts to an automorphism of $\widetilde{S}$ preserving the strict transform of the divisor $B$ and the $h$-exceptional locus. Write the group homomorphism as
$$\lambda:\Aut(S,B)\to\Aut(\tilde{S},\Delta_N)$$
Set ${m_0}:=m_0'N$. Then by the standard $\Aut(S,B)$-invariant inclusion
\[
\begin{aligned}
V_{mm_0}:=H^0(S,mm_0'N(K_S+B))
&=H^0(\widetilde{S},mm_0'N(K_{\widetilde{S}}+D_N)) \\
&\subseteq H^0(\widetilde{S},mm_0'NK_{\widetilde{S}}+(mm_0'N-mm_0')\lceil D_N\rceil) \\
&=H^0(\widetilde{S},mm_0'N(K_{\widetilde{S}}+\Delta_{N})),
\end{aligned}
\]
we know that the $A$-pluricanonical representation of $(S,B)$
$$\rho_{mm_0}:\Aut(S,B)\to\Aut(V_{mm_0})$$
is the composition of $\lambda$ and the $A$-pluricanonical representation for $(\tilde{S},\Delta_N)$ restricted on the subspace $V_{mm_0}$, and therefore is finite for every integer $m$.

Therefore, by replacing $(S,B)$ with $(\tilde{S},\Delta_N)$, we may assume $B$ has standard coefficients. By applying Lemma \ref{emmp}, we may assume  $K_S+B$ is nef and thus semi-ample by abundance for surfaces (see \cite[Theorem 1.2]{abundanceforlcsurf}).

To conclude, we need to prove the finiteness of $A$-pluricanonical representations for a two-dimensional Kawamata log terminal pair $(S,B)$ with standard coefficients and $K_S+B$ semi-ample. By Theorem \ref{main thm of kappa one} and \ref{main thm of cy surf}, it remains to prove for $(S,B)$ with $\kappa(K_S+B)=2$. In this case, the group $\Aut(S,B)$ is finite by \cite[Corollary 9.5]{Pat17}. This completes the proof.
\end{proof}

The main theorem of this section follows from the theorem below.

\begin{thm}\label{main thm: finiteness of pluri for peff lc surf}
Let $(S,B)$ be a projective log canonical pair of dimension $2$ over an algebraically closed field $k$ of characteristic $p >42$. Assume that $K_S+B$ is nef. Then there exists a positive integer $m_0$ such that the image of the $B$-pluricanonical representation
    \[
        \rho_{mm_0}:\Bir(S,B)\longrightarrow
        \Aut\big( H^0(S,mm_0(K_S+B))\big)
    \]
is finite for every positive integer $m$.
\end{thm}

\begin{proof}
We may assume $S$ is irreducible. By the abundance for surfaces (see \cite[Theorem 1.2]{abundanceforlcsurf}), we know that $K_S+B$ is semi-ample. Let
\[
f\colon S\to R
\]
be the morphism induced by $K_S+B$ such that $f_*\Oo_S=\Oo_R$. In fact, $f$ is induced by the Stein factorization of the morphism $S\to \proj\oplus_{m\geq0}H^0(S,mm_0(K_S+B))$ for a sufficiently divisible integer $m_0>0$.

Write $m_0(K_S+B)=f^*(H)$, where $H$ is an ample divisor on $R$. Notice that $\operatorname{Bir}(S,B)$ acts on $R$ as automorphisms given the fact that $g_*(m_0(K_S+B))=m_0(K_S+B)$ for every element $g\in \Bir(S,B)$. Denote the group morphism as $l:\Bir(S,B)\to \Aut(R)$.

\noindent
\textbf{Case 1: $\kappa(S,K_S+B)=2$.}

Then $f$ is a birational morphism. We put $\Xi:=f_*B$. Then
\[
K_S+B=f^*(K_R+\Xi)
\]
and $\Xi$ is $\operatorname{Bir}(S,B)$-invariant. Thus, we may reduce to the finiteness of $A$-pluricanonical representations of $(R,\Xi)$, which is a log canonical pair. Since $K_R+\Xi$ is ample, we conclude by \cite[Corollary 9.5]{Pat17}.

\medskip
\noindent
\textbf{Case 2: $\kappa(S,K_S+B)=1$.}
By \cite[Proposition 5.9]{tanaka2015abundancetheoremsemilog}, there exists a proper birational morphism $g\colon S'\longrightarrow S$ from a smooth surface \(S'\) such that \((S',B')\) is divisorial log terminal and $K_{S'}+B'=g^*(K_S+B)$. By Remark \ref{Main Thm equi under B-bir}, we may replace \((S,B)\) by \((S',B')\).

We first suppose that the horizontal part of of the divisor $\lfloor B\rfloor$ with respect to $f$ is non-zero, that is, $\lfloor B\rfloor ^h\neq 0$. Let \(D\) be an irreducible component of $\lfloor B\rfloor^h$. After possibly blowing up, we may assume $D$ is smooth. Define \(B_D\) by adjunction:
\[
K_D+B_D
 :=
 (K_S+D+(B-D))|_D
 =
 K_D+\operatorname{Diff}_D(B-D).
\]
Then \((D,B_D)\) is log canonical by \cite[Proposition~5.7]{tanaka2015abundancetheoremsemilog}. Note that we may assume the group $\Bir(S,B)$ preserves every component of $B^h$. Since if we have a $B$-birational map
\[\begin{tikzcd}
	& T & \\
	S && S \\
	R && R
	\arrow["\alpha"', from=1-2, to=2-1]
	\arrow["\beta", from=1-2, to=2-3]
	\arrow["g", dashed, from=2-1, to=2-3]
	\arrow["f"', from=2-1, to=3-1]
	\arrow["f", from=2-3, to=3-3]
	\arrow["{l(g)}", from=3-1, to=3-3]
\end{tikzcd}\]
then after considering a common log resolution, we have $\alpha^*D$ dominate $R$ and thus so does $\beta_*(\alpha^*D)$. After passing to the restricted representations on the subgroup $\Bir(S,B)_F$, we may assume every element in $\Bir(S,B)$ fixes every component of $B^h$. To see this, we consider the action of $\Bir(S,B)$ on the finitely many components of $B^h$. Then we know that $\Bir(S,B)_F$ is a normal subgroup of $\Bir(S,B)$ having only finitely many cosets.

Since \(D\to R\) is surjective and \(m_0(K_S+B)\sim f^*H\), restriction to \(D\) induces an injective map
\[
H^0\bigl(S,\mathcal O_S(m_1m_0(K_S+B))\bigr)
\hookrightarrow
H^0\bigl(
 D,\mathcal O_D(m_1m_0(K_D+B_D))
\bigr).
\]
Consequently, the assertion follows from the one-dimensional case, see Theorem \ref{1-dim subklt fin rep}.

We may therefore assume that $ \lf B\rfloor ^h=0$. Then $\lfloor B\rfloor$ is vertical. Define
\[
\Sigma:=\sum_{q\in f(\lfloor B\rfloor)} a_qf^*q
\]
with $$a_q=\max\{\frac{1}{\mult_Cf^*(q)}|C\subset \lfloor B\rfloor\text{ is an irreducible divisor with } f(C)=q\}.$$
Define a divisor $A:=\sum_{q\in f(\lfloor B\rfloor)} a_q q$. Then \(\Sigma=f^*A\). After passing to a normal subgroup of $\Bir(S,B)$ with finite index in $\Bir(S,B)$, we may assume that \(\Sigma\) is a \(\operatorname{Bir}(S,B)\)-invariant Weil divisor. In fact, $\Bir(S,B)$ preserves log canonical centers and we may take the subgroup consisting of elements $g\in \Bir(S,B)$ with $l(g)$ fixing $f(\lfloor B\rfloor)$, which has finite index in $\Bir(S,B)$. Note that the divisor $B':=\Sigma-\lfloor B\rfloor$ is effective and does not contain any whole fiber. By \cite[Lemma]{Bombieri_Mumford_1977}, the intersection matrix of each connected component of $B'$ is negative definite if $B'\neq 0$.

For any rational number $\epsilon>0$ that is sufficiently small so that $\frac{1}{m_0}H-\epsilon A$ has positive degree, we run the $(K_S+B-\epsilon\lfloor B\rfloor)$-MMP on $S$ over $R$. Note that $(S,B-\epsilon \lfloor B\rfloor)$ is Kawamata log terminal. Given $K_S+B-\epsilon\lfloor B\rfloor=(K_S+B-\epsilon \Sigma)+\epsilon B'$ with $K_S+B-\epsilon \Sigma$ being $f$-trivial and the intersection matrix of each connected component of $B'$ being negative definite, we can conclude that the Zariski decomposition for $K_S+B-\epsilon\lfloor B\rfloor$ has $K_S+B-\epsilon \Sigma$ as the positive part and $\epsilon B'$ as the negative part. Therefore, from the proof of Lemma \ref{emmp}, if we denote the result of the $(K_S+B-\epsilon\lfloor B\rfloor)$-MMP as a birational morphism $u:S\to T$ with $K_T+B_T:=u_*(K_S+B-\epsilon\Sigma)$, then $u^*(K_T+B_T)=K_S+B-\epsilon \Sigma$. Apply Theorem \ref{main thm: finite pluri repn for Kawamata log terminal} and Remark \ref{Main Thm equi under B-bir}, we have the finiteness of the image of the following representation

$$\rho_{m,\epsilon}: \Bir(S,B-\epsilon \Sigma)\to\Aut\big( H^0(m(K_S+B-\epsilon\Sigma))\big)$$
for any sufficiently divisible integer $m>0$ (we denote by $\rho_{m,\epsilon}$ to treat different $\epsilon$). Since the divisor $\Sigma$ is $\Bir(S,B)$-invariant, we have a natural inclusion
$$\Bir(S,B)\subset \Bir(S,B-\epsilon \Sigma)$$
and $\rho_{m,\epsilon}(\Bir(S,B))$ is finite. We fix $\epsilon_0>0$ that is sufficiently small and take $m_0>0$ such that $\rho_{mm_0,\epsilon_0}(\Bir(S,B)), \rho_{mm_0,\frac{\epsilon_0 }{2}}(\Bir(S,B))$ are finite for any positive integer $m$. Therefore we may take a $\Bir(S,B)$-invariant non-zero section $s\in H^0(m_0N(K_S+B-\epsilon_0\Sigma))$ for some integer $N$. In fact, if we denote $\rho_{m_0,\epsilon_0}(\Bir(S,B))=\{g_i\}_{i=1}^{N}$, then
$$s:=\prod_{i=1}^N (g_it)\in H^0(m_0N(K_S+B-\epsilon_0\Sigma))$$
is $\Bir(S,B)$-invariant for any non-zero section $t\in H^0(m_0(K_S+B-\epsilon_0\Sigma))$. The multiplication by $s^m$ induces a $\Bir(S,B)$-invariant inclusion
\[
\begin{aligned}
    H^0(S, mm_0N(K_S+B))\overset{\cdot s^{m}}{\hookrightarrow} & H^0(mm_0N(K_S+B)+mm_0N(K_S+B-\epsilon_0\Sigma))\\
    & =H^0\bigl(2mm_0N(K_S+B-\frac{\epsilon_0}{2}\Sigma)\bigl)
\end{aligned}
\]
Then the finiteness of $m$-th $B$-pluricanonical representation for $(S,B)$ follows from the finiteness of $\rho_{2mm_0N,\frac{\epsilon_0 }{2}}(\Bir(S,B))$.

\medskip
\noindent
\textbf{Case 3: $\kappa(S,K_S+B)=0$.}
If $(S,B)$ is Kawamata log terminal, we may conclude by Theorem \ref{main thm: finite pluri repn for Kawamata log terminal}. If not, by \cite[Proposition 5.9]{tanaka2015abundancetheoremsemilog}, we may replace $(S,B)$ with its divisorial log terminal model. In this case, we have $\lfloor B\rfloor\neq0$. Now take a non-zero preadmissible section $s\in \operatorname{PA}(S,mm_0(K_S+B))$ (see Definition \ref{def:admissible-section} for the definition) for any positive integer $m$ and sufficiently divisible positive integer $m_0$, whose existence is guaranteed by Claim $(A_2)$ in the proof of Theorem \ref{cor:slc-threefold-abundance} (note that Theorem \ref{main thm: finiteness of pluri for peff lc surf & Cartier ind} is not needed  to conclude Claim $(A_2)$). We will show that $g_*s=s$ for any $g\in\Bir(S,B)$, from which we know that the $B$-pluricanonical representations are trivial in this case. It suffices to show \(g_*s|_{\llcorner B \lrcorner} = s|_{\llcorner B \lrcorner}\). This follows from Lemma \ref{bir stable SG}.
\end{proof}

Now we are ready to prove the main theorem of this section.

\begin{proof}[Proof of Theorem \ref{main thm: finiteness of pluri for peff lc surf & Cartier ind}]
    We may assume $S$ is irreducible. Fix any positive integer $m$ such that $m(K_S+B)$ is Cartier. By Theorem \ref{main thm: finiteness of pluri for peff lc surf}, there exists a positive integer $m_0$, such that the following representation
\[
        \rho_{mm_0}:\Bir(S,B)\longrightarrow
        \Aut\big( H^0(S,mm_0(K_S+B))\big)
\]
is finite. After replacing $\Bir(S,B)$ with the kernel of $\rho_{mm_0}$, which has finite index, we may assume that $\Bir(S,B)$ acts on $ H^0(S,mm_0(K_S+B))$ trivially.

For any $s\in H^0(S,m(K_S+B))$, we have $s^{m_0}\in H^0(S,mm_0(K_S+B))$. Thus, for any $B$-birational map $g\in\Bir(S,B)$, we know that $g_*(s)^{m_0}=g_*(s^{m_0})=s^{m_0}$, which implies that $g_*(s)=as$ for some $a\in k^\times$ with $a^{m_0}=1$. Therefore, $\rho_m$ has finite image.
\end{proof}

\section{Reduced boundaries on divisorial log terminal threefolds}
Following Fujino’s strategy, we analyze the reduced boundary \(\lfloor B\rfloor\) of a divisorial log terminal pair in this section. This boundary records the gluing data arising from normalization and supplies the induction step. The main result of this section is Proposition \ref{reduced bdry of dlt 3 folds}, which corresponds to \cite[Proposition 2.1]{Fuj00}.

First, we need the following connectedness principle.

\begin{thm}[cf.{\cite[Theorem 2.15]{NT20}, \cite[Lemma 3.4]{Filipazzi2020ConnectednessPF}}]\label{connectedness for bir 3 dim char 5}
    Let $f : X \to V$ be a projective birational morphism of normal quasiprojective irreducible varieties of dimension at most $3$ over an algebraically closed field $k$ of characteristic $p>5$. Let $(X, B)$ be a sub-pair such that $-(K_X +B)$ is $f$-nef and $f_*B$ is effective. Then the induced morphism $\operatorname{Nklt}(X, B)\to V$ has connected fibers. 
\end{thm}

The following connectedness principle for threefolds was proved in \cite{Filipazzi2020ConnectednessPF}, which will play an important role in lifting sections (see Lemma \ref{lem:extension-from-boundary}), and we state it as follows.

\begin{thm}[cf.{\cite[Proposition 3.7]{Filipazzi2020ConnectednessPF}}]\label{connectedness for contraction 3 dim char 5}
    Let $f : X \to S$ be a projective morphism between quasi-projective irreducible varieties defined over an algebraically closed field $k$ of characteristic $p>5$. Assume that $\dim X \leq 3$ and that $(X, B)$ is a pair with rational coefficients and $-(K_X + B)$ nef over $S$. Assume that $\operatorname{Nklt}(X, B)$ is vertical over $S$ and that $f^{-1}(s) $ is connected as a $k(s)$-scheme. Then $f^{-1} (s) \cap \operatorname{Nklt}(X, B)$ is connected as a $k(s)$-scheme.
\end{thm}

First, we prove the following lemma, which investigates the case of Mori-fiber space in Proposition \ref{reduced bdry of dlt 3 folds}.

\begin{lemma}\label{reduced bdry for mori fiber sp.}
    Let $(Z,B)$ be a $\Qq$-factorial log canonical n-fold over an algebraically closed field $k$ of characteristic $p>42$ with $n=2$ or $3$. Let $h:Z\to R$ be a projective surjective morphism onto a normal variety $R$ such that $h_*\Oo_Z=\Oo_R$. Assume the following conditions:
\begin{enumerate}[labelsep=0.5em]
    \item  $\lfloor B\rfloor\neq0$ and \((Z, B - \epsilon \lfloor B \rfloor)\) is Kawamata log terminal, where \(\epsilon\) is a small positive rational number;
    \item \(K_Z + B \sim_{\Qq,h} 0\);
    \item there is a \((K_Z + B - \epsilon \lfloor B \rfloor)\)-extremal Fano contraction \(u : Z \rightarrow V\) over \(R\) such that \(\dim V = n - 1\).
\end{enumerate}

Then the horizontal part \(\lfloor B\rfloor^h\) of \(\lfloor B \rfloor\) with respect to the morphism $u$ is one of the following:

\begin{enumerate}[label=(\alph*), labelsep=0.5em]
    \item \(\lfloor B\rfloor^h = D_1\), which is irreducible, and \(\deg[D_1 : V] = 2\);
    \item \(\lfloor B\rfloor^h = D_1\), which is irreducible, and \(\deg[D_1 : V] = 1\);
    \item \(\lfloor B\rfloor^h = D_1 + D_2\), such that \(D_i\) is irreducible and \(\deg[D_i : V] = 1\) for \(i = 1, 2\).
\end{enumerate}
    
In the cases (a) and (c), the number of connected components of \(\lfloor B \rfloor \cap h^{-1}(r)\) is at most two for every \(r \in R\).

In the case (b), \(\lfloor B \rfloor \cap h^{-1}(r)\) is connected for every \(r \in R\).

In the case (a), there is a B-birational map \(\iota : (D_1, \operatorname{Diff}(B - D_1)) \dashrightarrow (D_1, \operatorname{Diff}(B - D_1))\) over \(V\) such that \(\iota^2 = \operatorname{id}\).

In the case (c), there is a $B$-birational map $\iota:(D_1,\Diff(B-D_1))\dashrightarrow(D_2,\Diff(B-D_2))$ over $V$.
\end{lemma}

\begin{proof}
Since \(\lfloor B \rfloor\) is \(u\)-ample by the assumptions (2) and (3), we have \(B^h \neq 0\). We can prove as the proof of Claim \ref{claim:kernel} to conclude the geometric normality of the generic fiber of $u$. We restate the proof as follows.  Take the generic fiber $F$, which is one-dimensional, normal and thus Gorenstein by the fact that $Z$ is normal. Suppose $F$ is not geometrically normal. We consider the normalization of $F_0$, the reduced subscheme of $\bar F:=F\otimes _{\spec k(V)}\overline{k(V)}$, and denote it by $\pi: T\to F_0$. Denote the composition of the morphisms $T\xrightarrow{\pi} F_0\hookrightarrow \bar F\to F$ as $\phi$. By \cite[Theorem 1.1]{PW17}, we have
$$K_T+(p-1)D\sim \phi^*K_{F}.$$
Take the degree and we have
$$-2\leq\deg K_T\leq -\deg (p-1)D\leq -p+1,$$
which is impossible since $p>42$.

Geometric normality is an open condition (cf.\cite[II, 6.9.1; III, 12.2.4(iv)]{EGAIV}). So the general fibers of \(Z \rightarrow V\) are geometrically normal and thus smooth by the fact that $k$ is algebraically closed. Also, we have the general fibers of $u$ are \(\mathbb{P}^1\) since $u$ is an extremal Fano contraction. This implies that \(\deg[\lfloor B\rfloor^h : V] \leq 2\) given that \(K_Z + B\) is numerically \(h\)-trivial. Thus, one of $(a)$, $(b)$, $(c)$ holds.

If we are in the case $(a)$ or $(c)$, then the horizontal part of $B$ has reduced coefficients, that is, $\lfloor B\rfloor^h=B^h$. Suppose first that we are in the case $(a)$, then the field extension $K(V)\to K(D_1)$ induced by $u|_{D_1}: D_1\to V$ has degree $2$ and is separable. Take an involution of $K(D_1)$ over $K(V)$ given by Galois theory. This involution induces a birational self-map $\iota:D_1\dashrightarrow D_1$ over $V$. If we are in the case $(c)$, the restricted morphisms $u|_{D_i}:D_i\to V$ are birational and $\iota$ is given by $\iota:= u|_{D_2}^{-1}\circ u|_{D_1}: D_1\dashrightarrow D_2$. It remains to prove that in either of these cases, $\iota$ is actually a $B$-birational (self-)map between $(D_i,\Diff_{D_i}(B-D_i))$ and $(D_j,\Diff_{D_j}(B-D_j))$. Here, we define $i=j=1$ if we are in case $(a)$, and $i=1,\ j=2$ if we are in case $(c)$.

When $n=2$, $\iota$ is an isomorphism. We remain to show that
\begin{equation}
\iota^*(\Diff_{D_j}(B-D_j))
=
\Diff_{D_i}(B-D_i).
\tag{$\star$}
\end{equation}
In case $(c)$, the proof follows from \cite[12.3.4 Theorem]{Kollr1992FlipsAA}. Here we give the proof for completeness. By taking a log resolution and contracting $(-1)$-curves in the fibers, we may assume $Z$ is smooth and $u$ is a $\mathbb{P}^1$-bundle. Now $(\star)$ follows from that $B$ consists of two sections and some fibers (with some coefficients). In case $(a)$, the ramification locus of $D_1\to V$ is $\iota$-invariant. On the other hand, over the maximal locus $U\subset D_1$ where the restricted morphism $U\to V$ is an unramified double cover, we verify $(\star)$ over the \'etale neighborhood $U\to V$. Namely, we consider the following base change.
\begin{center}
    \begin{tikzcd}
Z' \arrow[d, "u'"'] \arrow[r, "f"] & Z \arrow[d, "u"] \\
U \arrow[r, "u|_U"]                & V               
\end{tikzcd}
\end{center}
Then $f^*D_1=D_1'+D_2'$ with $D_1'$, $D_2'$ two sections over $U$. Define $B':=f^*B$, then $K_{Z'}+B'\sim_{\Qq, u'} 0$. Denote $\iota'$ to be the isomorphism between $D_1'$ and $D_2'$. Applying the proof of case $(c)$, we know that ${\iota'}^*(\Diff_{D_j'}(B'-D_j'))=\Diff_{D_i'}(B'-D_i')$ and that the equality $(\star)$ holds by pushing forward.

When $n=3$, we take a common log resolution as follows:
\begin{center}
    \begin{tikzcd}
                                           & W \arrow[ld, "\alpha"'] \arrow[rd, "\beta"] &                \\
D_i \arrow[rd] \arrow[rr, "\iota", dashed] &                                             & D_j \arrow[ld] \\
                                           & V                                           &               
\end{tikzcd}
\end{center}
We need to verify the identity $\alpha^*(K_{D_i}+\Diff_{D_i}(B-D_i))=\beta^*(K_{D_j}+\Diff_{D_j}(B-D_j))$. By taking a vertical general hyperplane $H$ on $V$, we denote the inverse image of $D_i\cap f^{-1}H$ with respect to the log resolution to be $\tilde{W}$. Consider the restricted diagram as follows:
\begin{center}
\begin{tikzcd}
                                                         & \tilde{W} \arrow[ld, "\alpha|_{\tilde{W}}"'] \arrow[rd, "\beta|_{\tilde{W}}"] &                           \\
D_i\cap {f^{-1}H} \arrow[rd] \arrow[rr, "\iota|_{D_i\cap{f^{-1}H}}", dashed] &                                                      & D_j\cap{f^{-1}H} \arrow[ld] \\
                                                         & H                                                    &                          
\end{tikzcd}
\end{center}
Denote the Stein factorization of $W\to V$ as $W\xrightarrow{\pi} S\to V$, where $\pi$ is a birational morphism. By the $n=2$ case, we conclude that the difference
$$\alpha^*(K_{D_i}+\Diff_{D_i}(B-D_i))-\beta^*(K_{D_j}+\Diff_{D_j}(B-D_j))$$
supports on divisors that are exceptional over $S$. By negativity lemma, we know that the difference is zero since it is numerically trivial over $S$.
\end{proof}

We can prove the main result in this section.

\begin{proposition}\label{reduced bdry of dlt 3 folds}
    Let (X,B) be an irreducible $\Qq$-factorial divisorial log terminal pair over an algebraically closed field $k$ of characteristic $p>42$ with $\dim X\leq3$. Let $f:X\to R$ be a projective morphism onto a normal variety $R$ such that $f_*\Oo_X=\Oo_R$. Assume that $K_X+B$ is $\Qq$-linearly $f$-trivial and $\lfloor B\rfloor\neq 0$. Then one of the following holds.

\begin{enumerate}[labelsep=0.5em]
    \item $\operatorname{dim} R =0$.
    \begin{enumerate}[labelsep=0.5em]
        \item  $\lfloor B\rfloor$ is connected.
        \item  $\lfloor B\rfloor$ has two connected components $B_1$ and $B_2$, and there exists a rational map $v:X\dashrightarrow V$ onto a $\Qq$-factorial log canonical $(\dim X-1)$-fold $V$ with general fiber $\mathbb{P}^1$. Furthermore, there exists an irreducible component $D_i'\subset B_i$ such that there exists a $B$-birational map over $V$
        $$\iota:(D_1',\Diff(B-D_1'))\dashrightarrow(D_2',\Diff(B-D_2')).$$
    \end{enumerate}  
    
    \item $\operatorname{dim} R \geq1$.
    \begin{enumerate}[labelsep=0.5em]
        \item  $\lfloor B\rfloor\cap f^{-1}(r)$ is connected for every $r\in R$.
        \item  The number of connected components of $\lfloor B\rfloor\cap f^{-1}(r)$ is at most two for every $r\in R$. There exists a rational map $v:X\dashrightarrow V$ onto a $\Qq$-factorial log canonical (n-1)-fold $V$ with general fiber $\mathbb{P}^1$. The horizontal part $\lfloor B\rfloor^h$ of $\lfloor B\rfloor$, with respect to the map $v$, is one of the following:
        \begin{enumerate}[labelsep=0.5em]
            \item $\lfloor B\rfloor^h=D_1'$, which is irreducible, and the mapping degree $\deg[D_1':V]=2$; there is also a $B$-birational involution on $(D_1',\Diff(B-D_1'))$ over $V$;
            \item $\lfloor B\rfloor^h=D_1'+D_2'$ such that $D_i'$ is irreducible and there is a $B$-birational map
            $$(D_1',\Diff(B-D_1'))\dashrightarrow(D_2',\Diff(B-D_2')),$$
            which we also denoted as $\iota$.
        \end{enumerate}
    \end{enumerate}
\end{enumerate}
\end{proposition}

\begin{proof}
    If $f$ is birational, by Theorem \ref{connectedness for bir 3 dim char 5}, we are in the case $(1.a)$ or $(2.a)$. Thus, we may assume that $\dim R<\dim X$.
    
    We run a $(K_X+B-\epsilon\lfloor B\rfloor)$-MMP on $X$ over $R$ for $0<\epsilon\ll1$. The end result is a birational map $p:X\dashrightarrow Z$ over $R$. Let $h:Z\to R$ be the induced morphism. Since $K_X+B$ is numerically trivial over $R$, so is the divisor $K_Z+p_*B$ on $Z$. Also, $(Z,p_*B)$ is a $\Qq$-factorial log canonical pair by Lemma \ref{decsend CY under MMP} and $(Z,p_*B-\epsilon\lfloor p_*B\rfloor)$ is Kawamata log terminal.
    
    \textit{Step 1.} If $(Z,p_*B-\epsilon\lfloor p_*B\rfloor)$ is a  minimal model, then $K_Z+p_*B-\epsilon\lfloor p_*B\rfloor$ is $h$-nef and thus $-\lfloor p_*B\rfloor$ is $h$-nef since $K_Z+p_*B$ is $h$-trivial. If $\dim R=0$, then $\lfloor p_*B\rfloor=0$. By Lemma \ref{invariance of number of connected comp of reduced bdry under MMP}, we have $\lfloor B\rfloor=0$. This case is excluded. If $0<\dim R< \dim X$, then $\lfloor p_*B\rfloor\cap h^{-1}(r)$ is connected for every $r\in R$. In fact, if $\lfloor p_*B\rfloor\cap h^{-1}(r)$ is nonempty for some $r\in R$, then it must support on the whole fiber $h^{-1}(r)$, otherwise there would be an irreducible curve in the fiber intersecting with $\lfloor p_*B\rfloor\cap h^{-1}(r)$ but not contained in the support of $\lfloor p_*B\rfloor\cap h^{-1}(r)$, contradicting the fact that $-\lfloor p_*B\rfloor$ is $h$-nef. By Lemma \ref{invariance of number of connected comp of reduced bdry under MMP} we have the case $(2.a)$.

    \textit{Step 2.} If there exists an extremal Fano contraction $u:Z\to V$ over $R$, which is  $-(K_Z+p_*B-\epsilon\lfloor p_*B\rfloor)$ is $u$-ample. Let $g:V\to R$ be the induced morphism. We note that $h=g\circ u$.

    First, assume $\dim R=0$. If $\lfloor p_*B\rfloor$ is connected, we obtain case $(1.a)$ by Lemma \ref{invariance of number of connected comp of reduced bdry under MMP}. Thus, we can assume that $\lfloor p_*B\rfloor$ is not connected. Since $\lfloor p_*B\rfloor$ is $u$-ample, we have $\dim V=\dim X-1$, since an effective ample Weil divisor on a surface has connected support. Then by Lemma \ref{reduced bdry for mori fiber sp.}, we have the case (c) in Lemma \ref{reduced bdry for mori fiber sp.} and $u|_{D_i}\to V$ is finite, where $D_i$ are as in Lemma \ref{reduced bdry for mori fiber sp.}(c) for $i=1,2$. By Lemma \ref{invariance of number of connected comp of reduced bdry under MMP}, $\lfloor B\rfloor$ has two connected components $B_1$ and $B_2$. By Lemma \ref{decsend CY under MMP} and \ref{reduced bdry for mori fiber sp.}, we have the case $(1.b)$ by taking $D_i'$ to be the strict transform of $D_i$.

    Thus, we may assume $\dim R\geq1$.

    If $\dim V=\dim X-1$, similarly, by Lemmas \ref{decsend CY under MMP}, \ref{reduced bdry for mori fiber sp.}, and \ref{invariance of number of connected comp of reduced bdry under MMP}, we have that the horizontal part $\lfloor B\rfloor ^h=D_1'$ or $\lfloor B\rfloor ^h=D_1'+D_2'$ (rather than $\lfloor B\rfloor$ contains $D_i'$ as in the case $(1.b)$) and we get the case $(2.a)$ if $\lfloor B\rfloor ^h=D_1'$ with $\deg[D_1':V]=1$ and $(2.b)$ otherwise.

    If $\dim R=1$, $\dim V=1$, and $\dim X=3$, then $V\simeq R$. Since $u$ is extremal, and $\rho(Z/R)=1$, every vertical irreducible component of $\lfloor p_*B \rfloor$ is a pullback of a $\Qq$-Cartier $\Qq$-divisor on $R$, and every horizontal irreducible component of $\lfloor p_*B\rfloor$ is $h$-ample. Since an effective ample divisor on a surface must have connected support, $\lfloor p_*B\rfloor\cap h^{-1}(r)$ is connected for every $r\in R$. Thus, we have the case $(2.a)$.
\end{proof}

The following lemma is used in the proof of Proposition \ref{reduced bdry of dlt 3 folds}.

\begin{lem}\label{invariance of number of connected comp of reduced bdry under MMP}
    Given morphisms $f:X\to R$, $h:Z\to R$, $p:X\dashrightarrow Z$ as in the proof of Proposition \ref{reduced bdry of dlt 3 folds}, the number of connected components of $\lfloor B\rfloor\cap f^{-1}(r)$ is equal to the number of connected components of  $\lfloor p_*B\rfloor\cap h^{-1}(r)$ for every $r\in R$.
\end{lem}
\begin{proof}
    The birational map $p: X\dashrightarrow X_1\dashrightarrow X_2\dashrightarrow \cdots\dashrightarrow Z$ is a composition of flips and divisorial contractions. Note that $\lfloor B_i\rfloor$ is relatively ample for each flipping or divisorial contraction and $K_{X_i} + B_i$ is numerically trivial over $R$, where $B_i$ is the boundary on $X_i$ induced by $(X,B)$ and $f_i:X_i\to R$ the induced morphism. We prove a similar assertion for $X_i\dashrightarrow X_{i+1}$ step by step, and the conclusion follows. Note that the pairs \((X_i,B_i-\epsilon\lfloor B_i\rfloor))\) are klt by the construction, and thus for $\epsilon$ sufficiently small, we have \(\operatorname{Nklt}(X_i,B_i)=\operatorname{Supp}\lfloor B_i\rfloor\),

    Say one of the steps, $X_i\to X_{i+1}$, is a divisorial contraction, using Theorem \ref{connectedness for bir 3 dim char 5}, we know that the number of connected components of $\lfloor B_i\rfloor\cap f_i^{-1}(r)$ is equal to the number of connected components of $\lfloor B_{i+1}\rfloor\cap f_{i+1}^{-1}(r)$.

    If $X_i\to Y_{i+1}\xleftarrow{} X_{i+1}$ is a flip, we apply Theorem \ref{connectedness for bir 3 dim char 5} for the morphisms $X_i\to Y_{i+1}$ over $R$ and $X_{i+1}\to Y_{i+1}$ over $R$, and we know that the same conclusion holds in this step.
\end{proof}
\section{The abundance theorem for semi log canonical threefolds}

In this section, we prove the abundance conjecture for semi log canonical threefolds. We first recall the following definitions.

\begin{definition}\label{def : sdlt slc}
    Let $X$ be a reduced $S_2$ scheme.  We assume that it is pure $n$-dimensional and that it is regular or normal crossing in codimension 1. Let $B$ be an effective $\Qq$-divisor on $X$ such that $K_X+B$ is $\Qq$-Cartier.

    Let $X=\cup_iX_i$ be the decomposition into irreducible components, and let
$$\mu : X' := \coprod_iX'_i \rightarrow X = \cup X_i$$
be the normalization. A $\Qq$-divisor $B_i$ on $X_i'$ is defined by 
$$K_{X_i'}+B_i=\mu^*(K_X+B)|_{X_i'}$$

    We say that $(X,B)$ is \emph{semi log canonical} if $(X_i',B_i)$ is log canonical for all $i$.

    We say that $(X,B)$ is \emph{semi divisorial log terminal} if $X_i$ is normal and $(X'_i,B_i)$ is divisorial log terminal for all $i$.
\end{definition}
\begin{rmk}\label{sdlt of reduced part of dlt 3-folds}
    Let $(X,B)$ be a divisorial log terminal pair. In characteristic zero, the adjunction $(\lfloor B\rfloor,\Diff(B-\lfloor B\rfloor))$ is semi divisorial log terminal by the Grauert-Riemenschneider vanishing (see \cite[Corollary 5.52]{Kollár_Mori_1998} and \cite[Corollary 17.5]{Kollr1992FlipsAA}).
    
    If $X$ is a surface over an algebraically closed field of characteristic $p>0$, then adjunction $(\lfloor B\rfloor,\Diff(B-\lfloor B\rfloor))$ is semi divisorial log terminal by \cite[Corollary 5.8]{tanaka2015abundancetheoremsemilog}.
    
    If $X$ is a threefold over an algebraically closed field of characteristic $p>5$, by \cite[Theorem 3, 19]{bernasconi2023vanishing}, we obtain that $\lfloor B\rfloor$, the reduced part of the boundary, is Cohen-Macaulay and each irreducible component of $\lfloor B\rfloor$ is normal. Thus, the pair $(\lfloor B\rfloor, \Diff(B-\lfloor B\rfloor))$ is semi divisorial log terminal by adjunction (see \cite[Proposition 2.2]{DH16}).
\end{rmk}

Following \cite{Fuj00}, we introduce the notion of \emph{admissible} and \emph{preadmissible sections} for the proof of the abundance conjecture for semi log canonical $n$-folds.

\begin{definition}\label{def:admissible-section}
Let $(X,B)$ be a semi divisorial log terminal pair of dimension $n\leq 3$ over an algebraically closed field $k$ of characteristic $p>5$. Let $X=\cup_iX_i$ be the decomposition into irreducible components, and let $\mu : X^\nu := \coprod_iX_i \rightarrow X = \cup_i X_i$ be the normalization. Note that $(X_i,B_i)$ is divisorial log terminal. We define \emph{admissible} and \emph{preadmissible sections} inductively on dimension:
\begin{itemize}
\item     We say $s\in H^0\bigl(X, m(K_X+B)\bigr)$ is \emph {preadmissible} if the restriction $s|_{\coprod_i\lfloor B_i\rfloor}$ is admissible;
    \item We say $s\in H^0\bigl(X,m(K_X+B)\bigr)$ is \emph{admissible} if it is preadmissible and $ g^*(s|_{X_j})=s|_{X_i}$ for every $B$-birational map $g\colon (X_i,B_i)\dashrightarrow (X_j,B_j)$;
\end{itemize}
We define linear subspaces of $H^0\bigl(X,m(K_X+B)\bigr)$ as follows:
\[
\operatorname{A}(X,m(K_X+B))
:=
\{s\mid s \text{ is admissible}\}
\]
and
\[
\operatorname{PA}(X,m(K_X+B))
:=
\{s\mid s \text{ is preadmissible}\}.
\]
The zero-dimensional convention is important. For \(P=\coprod_{j=1}^r\operatorname{Spec} k\), set
$$
\operatorname{PA}(P)=k^r, \qquad \operatorname{A}(P)=\{(c,\ldots,c):c\in k\}.
$$
\end{definition}

The following lemma encapsulates the idea of reducing abundance conjecture for semi log canonical varieties to abundance conjecture for divisorial log terminal varieties, that is, finding enough preadmissible sections on the normalization.

\begin{lemma}\label{lem:pa-descends-normalization}
Let $(X,B)$ be a projective semi log canonical pair of dimension $n\leq 3$ over an algebraically closed field $k$ of characteristic $p>5$. Let $\nu:(Y,B_Y)\to (X,B)$ be the normalization. Let $(Z,B_Z)$ be a divisorial log terminal model of $(Y,B_Y)$, that is, there is a birational morphism $\mu:Z\to Y$ such that $\mu^*(K_Y+B_Y)=K_Z+B_Z$ and $(Z,B_Z)$ is divisorial log terminal. Assume $s\in \operatorname{PA}(Z,m(K_Z+B_Z))$ for a sufficiently divisible integer $m$. Then $s=\mu^*\nu^*t$ for some $t\in H^0(X,m(K_X+B))$.
\end{lemma}

\begin{proof}
This follows from \cite[Theorem 3.7]{tanaka2015abundancetheoremsemilog}.
\end{proof}

The following lemma is used to descend sections along proper morphisms with connected fibers.

\begin{lem} \label{connectedness fiber can descend section}
Assume $S$, $T$ are two varieties over an algebraically closed field $k$ of characteristic $p>0$. Let $f: S\to T$ be a proper surjective morphism with connected fibers. Then for every positive integer $e$ that is sufficiently large, the following holds:

For any line bundle $H$ on $T$, and any section $s\in H^0(S,f^*H)$, there exists a section $u\in H^0(T,p^eH)$ such that $s^{p^e}=f^*u$.
\end{lem}
\begin{proof}
    Consider the Stein factorization $f:S\xrightarrow{c} T'\xrightarrow{g} T$. We know that every closed fiber of the finite map $g$ has one $k$-point. The closed diagonal $T'\to T'\times_TT'$ therefore contains every closed point, hence is surjective. Thus $g$ is universally injective, and consequently a universal homeomorphism.

Therefore, there exists a positive integer $e$ such that the $e$-th iterated Frobenius of $T'$ factors through $g$, that is, there exists a morphism $h:T\to T^{'}$ such that we have the following commutative diagram
\begin{center}
            \begin{tikzcd}
            F^e:\ \ T^{'} \arrow[r, "g"] & T \arrow[r, "h"] & T^{'}
            \end{tikzcd}
\end{center}

By construction of $T^{'}$, we can find a section $t\in H^0(T',g^*H)$
$$s=c^*t$$
Consider $t^{p^e}=F^{e*}t=g^*(h^*t)\in H^0(T^{'},p^eg^*(H))$. Thus, we have $s^{p^e}=f^*(h^*t)$. Note that for a larger integer $e'$, the $e'$-th iterated Frobenius map factors through the $e$-th Frobenius map, and we can conclude the same way.
\end{proof}

We use the next lemma in the proof of Proposition \ref{prop:preadmissible-extension} when the reduced boundary $\lfloor B\rfloor$ does not dominate $R$.

\begin{lemma}\label{lem:extension-from-boundary}
Let $(X,B)$ be an irreducible $\mathbb{Q}$-factorial divisorial log terminal threefold over an algebraically closed field $k$ of characteristic $p>0$ with $K_X+B$ nef, and $\lfloor B\rfloor\neq 0$. Assume that
\[
f:=\Phi_{|m_0(K_X+B)|}\colon X\to R
\]
is a proper surjective morphism onto a normal variety $R$ induced by the linear system $|m_0(K_X+B)|$ for a sufficiently large and divisible integer $m_0$ such that $f_*\Oo_X=\Oo_R$, and that
\[
f(\lfloor B\rfloor)\subsetneq R.
\]
If there exist sections
\[
\{s_i\}_{i=1}^{d}
\subset
H^0\bigl(\lfloor B\rfloor,m(K_X+B)|_{\lfloor B\rfloor}\bigr)
\]
without common zeros for an integer $m$ divisible by $m_0$, then there exist sections
\[
\{u_i\}_{i=1}^{l}
\subset
H^0\bigl(X,p^em(K_X+B)\bigr)
\]
for some integer $e$ such that
\[
u_i|_S=
\begin{cases}
s_i^{p^e} & \text{for } 1\leq i\leq d,\\
0 & \text{for } d<i\leq l,
\end{cases}
\]
and $\{u_i\}_{i=1}^{l}$ have no common zeros.
\end{lemma}

\begin{proof}
Define $S:=\lfloor B\rfloor$. There is an ample $\mathbb{Q}$-Cartier $\mathbb{Q}$-divisor $H$ on $R$ such that
\[
K_X+B\sim_{\mathbb{Q}} f^*H.
\]

By Theorem \ref{connectedness for contraction 3 dim char 5}, we have $S\cap f^{-1}(r)$ is connected for all $r\in R$. Thus, apply Lemma \ref{connectedness fiber can descend section} with the restricted morphism $f|_S:S\to f(S)$ and the divisor $H|_{f(S)}$, we can get sections $\{t_i\}_{i=1}^d\subset H^0\bigl(f(S),p^emH|_{f(S)}\bigr)$ such that
$$s_i^{p^e}=f|_S^*(t_i)$$
for all $i$.

Now we consider the following restriction
$$r:H^0(R,p^emH)\to H^0\bigl(f(S),p^emH|_{f(S)}\bigr).$$
Since $H$ is ample, we may assume the integer $e$ to be sufficiently large and divisible so that the restriction is surjective. Choose $v_i\in H^0(R,p^emH)$ such that $r(v_i)=t_i$. Then we can lift the section $s_i^{p^e}$ to
$$\{u_i:=f^*v_i\}_{i=1}^d\subset H^0\bigl(X,p^em(K_X+B)\bigr).$$
 Also, we may assume $e$ is sufficiently large such that
\[
\mathcal{O}_R(p^emH)\otimes \mathscr{I}_{f(S)}
\]
is generated by global sections, where $\mathscr{I}_{f(S)}$ is the ideal sheaf of $f(S)$. Therefore we can get sections $\{u_i\}_{i=1}^{l}$ with the required properties.
\end{proof}

The next proposition is the main part of the proof of the abundance for semi log canonical pairs. The formulation is slightly different from \cite[Proposition 4.5]{Fuj00} because of positive characteristic (see Lemma \ref{connectedness fiber can descend section}).

\begin{proposition}\label{prop:preadmissible-extension}
Let $(X,B)$ be a projective $\mathbb{Q}$-factorial divisorial log terminal pair over an algebraically closed field $k$ of characteristic $p>42$ with $\dim X\leq 3$. Let $m$ be a sufficiently large and divisible integer, such that $m\in 2\mathbb{Z}$. Assume that:
\begin{enumerate}
    \item $K_X+B$ is nef;
    \item$\operatorname{A}\bigl(
    \lfloor B\rfloor,
    m(K_X+B)|_{\lfloor B\rfloor}
    \bigr)$ generates $\mathcal{O}_{\lfloor B\rfloor}
    \bigl(
    m(K_{\lfloor B\rfloor}
    +\operatorname{Diff}(B-\lfloor B\rfloor))
    \bigr)$.
\end{enumerate}
Then there exists an integer $e\geq0$, such that every section
$$s^{p^e}\in \operatorname{A}\bigl(\lfloor B\rfloor,mp^e(K_X+B)|_{\lfloor B\rfloor}\bigl)$$
with $s\in \operatorname{A}\bigl(
\lfloor B\rfloor, m(K_X+B)|_{\lfloor B\rfloor}\bigl)$ can be lifted to $X$. Equivalently, there exists $u\in\operatorname{PA}\bigl(
X,mp^e(K_X+B)\bigl)$, such that $u|_{\lfloor B\rfloor}=s^{p^e}$.

Furthermore, we have $\operatorname{PA}(X,mp^e(K_X+B))$ generates $\mathcal{O}_X(mp^e(K_X+B))$.
\end{proposition}

\begin{proof}
We may assume that $X$ is irreducible. By abundance for log canonical threefolds (see \cite[Theorem 1.1]{bernasconi2026abundancethreefoldspositivecharacteristic}), we know that $K_X+B$ is semi-ample. Consider the fibration induced by the linear system $|m_0(K_X+B)|$
\[
f:=\Phi_{|m_0(K_X+B)|}\colon X\to R
\]
for a sufficiently large and divisible integer $m_0$. If $\lfloor B\rfloor=0$, then the proposition is immediate. Thus we may assume $\lfloor B\rfloor\neq 0$.

We have the following possibilities by Proposition~\ref{reduced bdry of dlt 3 folds}:
\begin{enumerate}[label=(\arabic*)]
    \item $\kappa(X,K_X+B)=0$ and $\lfloor B\rfloor$ is connected;
    \item $\kappa(X,K_X+B)\geq 1$ and $f^{-1}(r)\cap \lfloor B\rfloor$ is connected for every $r\in R$ and $f(\lfloor B\rfloor)=R$;
    \item $\kappa(X,K_X+B)\geq 1$ and $f^{-1}(r)\cap \lfloor B\rfloor$ is connected for every $r\in R$ and $f(\lfloor B\rfloor)\subsetneq R$;
    \item $f^{-1}(r)\cap \lfloor B\rfloor$ is not connected for some $r\in R$.
\end{enumerate}

\medskip
\noindent
\textbf{Case 1:}
Consider the exact sequence
\[
0
\longrightarrow
H^0\bigl(X,m(K_X+B)-\lfloor B\rfloor)\bigr)
\longrightarrow
H^0\bigl(X,m(K_X+B)\bigr)
\]
\[
\longrightarrow
H^0\bigl(
\lfloor B\rfloor,
m(K_X+B)|_{\lfloor B\rfloor}
\bigr)
\longrightarrow
\cdots .
\]
Since $\kappa(X,K_X+B)=0$, we have $H^0\bigl(X,m(K_X+B)-\lfloor B\rfloor\bigr)=0$, and the second and third terms are one-dimensional. Thus we obtain the result.

\medskip
\noindent
\textbf{Case 2:}
Consider the restricted morphism $f|_{\lf B\rf}:\lf B\rf\to R$. By Lemma \ref{connectedness fiber can descend section}, there exists a positive integer $e$, such that we can lift $s^{p^e}$ to a preadmissible section on $X$ for any section $s\in \operatorname{A}(\lfloor B\rfloor,m(K_X+B)|_{\lfloor B\rfloor})$. The remaining assertion follows the fact that $\lfloor B\rfloor$ dominates $R$.

\medskip
\noindent
\textbf{Case 3:}
By Lemma~\ref{lem:extension-from-boundary}, we can lift
\[
s^{p^e}\in
\operatorname{A}\bigl(
\lfloor B\rfloor,
p^em(K_X+B)|_{\lfloor B\rfloor}
\bigr)
\]
to a section
\[
u\in \operatorname{PA}(X,p^em(K_X+B))
\]
if $e$ is a sufficiently large and divisible integer. By Lemma \ref{lem:extension-from-boundary}, $\operatorname{PA}(X,p^em(K_X+B))$ generates $\mathcal{O}_X(p^em(K_X+B))$ if $e$ is a sufficiently large and divisible integer.

\medskip
\noindent
\textbf{Case 4:}
We run the log $(K_X+B-\epsilon\lfloor B\rfloor)$-MMP over $R$ as in the proof of Proposition~\ref{reduced bdry of dlt 3 folds}. 
Let
\[
u\colon (X',B')\to (V,P)
\]
be the last step of the MMP. By Proposition \ref{reduced bdry of dlt 3 folds}, $u$ is a Fano contraction to an $(n-1)$-dimensional variety $V$ in this case. Write $K_{X'}+B'=u^*A$ for a divisor $A$ on $V$.

Note that $(X^{'}, B^{'})$ is $\Qq$-factorial log canonical, with the divisor $K_{X^{'}}+B^{'}$ $\Qq$-linearly trivial over $R$.  Also, $(X, B)$ and $(X^{'}, B^{'})$ are $B$-birational by Lemma \ref{decsend CY under MMP}. We denote the map to be $p\colon X\dashrightarrow X'$.

Let
\[
\alpha\colon (Y,B_Y)\to (X,B)
\qquad\text{and}\qquad
\beta\colon (Y,B_Y)\to (X',B')
\]
be a common log resolution of the $B$-birational map $p\colon X\dashrightarrow X'$, such that
\[
K_Y+B_Y=\alpha^*(K_X+B)
\qquad\text{and}\qquad
K_Y+B_Y=\beta^*(K_{X'}+B').
\]

Consider the induced morphisms
$$\alpha|_{B_Y^{=1}}: {B_Y^{=1}}\to \lfloor B\rfloor \qquad\text{and}\qquad
\beta|_{B_Y^{=1}}: {B_Y^{=1}}\to \lfloor B^{'}\rfloor.$$
For every section $s\in \operatorname{A}\bigl(\lfloor B\rfloor, m(K_X+B)|_{\lfloor B\rfloor}\bigl)$, we first pull it back to a section $$s'':=\alpha|_{B_Y^{=1}}^*(s)\in H^0\bigl(B_Y^{=1}, m(K_Y+B_Y)|_{B_Y^{=1}}\bigl)$$
on $Y$. Since the morphisms $\alpha|_{B_Y^{=1}},\beta|_{B_Y^{=1}}$ have connected fibers by Theorem \ref{connectedness for bir 3 dim char 5}, we obtain an integer $e$ and a section
$$s'\in H^0\bigl(\lfloor B^{'}\rfloor,p^em(K_{X^{'}}+B^{'})|_{\lfloor B^{'}\rfloor}\bigl),$$
such that
$$s''^{p^e}=\beta|_{B_Y^{=1}}^*(s').$$
by applying Lemma \ref{connectedness fiber can descend section}. Note that the choice of the integer $e$ depends only on the morphisms $\alpha|_{B_Y^{=1}},\beta|_{B_Y^{=1}}$.

Write $\lfloor B'\rfloor=B'^{h}+B'^v$ as the decomposition of horizontal divisors and vertical divisors with respect to $u$. Denote the strict transform of $B'^h$ on $X$ by $B^h$. We claim that $p|_{B^h}$ induces a $B$-birational map between $( B^h, \Diff( B-B^h))$ and $( B'^h, \Diff( B'-B'^h))$. Again, the MMP $X\dashrightarrow X'$ is over $R$, and the morphism $B^h\to V$, $B'^h\to V$ are generically finite, implying that none of the components of $B^h$ is contracted by $p$, and that none of the components of $B'^h$ is extracted by $p$. By pulling back to a common log resolution and adjunction, we conclude the claim.

Thus we know that $s'|_{B'^h}$ is $B$-birational involution invariant over $V$. Therefore, it descends to a section $t\in H^0(V,p^emA)$. Now we take $w'=u^*t\in H^0(X',p^em(K_{X'}+B'))$, we will prove that $w'^{p^{e'}}$ is a lift of the section $s'^{p^{e'}}$, i.e. $w'^{p^{e'}}|_{\lfloor B'\rfloor}=s'^{p^{e'}}$ for some positive integer $e'$.

By construction, $w'^{p^{e'}}|_{B'^h}=s'^{p^{e'}}|_{B'^h}$. Now we need to find a positive integer $e'$ such that $w'^{p^{e'}}|_{B'^v}=s'^{p^{e'}}|_{B'^v}$. Since $u$ is extremal and $B'^v$ is vertical, we have $B'^v=u^*D$ for some divisor $D$ on $V$. Write $D=\sum a_iD_i$ and $E_i=u^{-1}D_i$ with the reduced scheme structure. It is enough to check $w'^{p^{e'}}|_{E_i}=s'^{p^{e'}}|_{E_i}$.

Let $W_i$ be a component of $B'^h\cap E_i$ such that $W_i$ dominates $D_i$. We can always take such $W_i$ since $u(B'^h)=V$. Apply Lemma \ref{connectedness fiber can descend section} to the morphism $u|_{E_i}:E_i\to D_i$, and we obtain an integer $e'$ as in the conclusion of Lemma \ref{connectedness fiber can descend section}. By the choice of $e'$, we have the following diagram, where the left vertical map is the descending map, and the right vertical map is the pullback.

\begin{center}

\begin{tikzcd}[row sep=3.5em, column sep=5em]
    (H^0(E_i, p^e m(K_{X'} + B')|_{E_i}))^{p^{e'}}
    \arrow[d, "\iota"]
        \arrow[r, "r"] &
    H^0(W_i,p^{e+e'} m(K_{X'} + B')|_{W_i})
         \\
    H^0(D_i, p^{e+e'} mA|_{D_i})
        \arrow[r, "id"] &
    H^0(D_i, p^{e+e'} mA|_{D_i})
        \arrow[u, "\gamma"]
\end{tikzcd}

\end{center}

The maps $\iota$, $\gamma$ are injective since $u|_{E_i}$, $u|_{W_i}$ are dominant. Thus the map $r$ is injective. Since $w'|_{B'^h}=s'|_{B'^h}$, we have $w'^{p^{e'}}|_{W_i}=s'^{p^{e'}}|_{W_i}$, which implies that the inverse images $w'^{p^{e'}}|_{E_i}=r^{-1}(w'^{p^{e'}}|_{W_i})$ and $s'^{p^{e'}}|_{E_i}=r^{-1}(s'^{p^{e'}}|_{W_i})$ are equal.

Under the isomorphism
$$ H^0(Y,p^{e+e'}m(K_Y+B_Y))\simeq H^0(X',p^{e+e'}m(K_{X'}+B')),$$
we obtain a section $w''=\beta^*(w'^{p^{e'}})\in  H^0(Y,p^{e+e'}m(K_Y+B_Y))$, whose restriction on $B_Y^{=1}$, $w''|_{B_Y^{=1}}$ is equal to $s''^{p^{e+e'}}=\beta|_{B_Y^{=1}}(s'^{p^{e'}})$. By the choice of the integer $e$, the section $s''^{p^{e+e'}}$ descends to a section on $\lfloor B\rfloor$, which can only be $s^{p^{e+e'}}$ by the dominance of the morphism $\alpha|_{B_Y^{=1}}:B_Y^{=1}\to \lfloor B\rfloor$. And we can also descend the section $w''$ to a section on $X$, denoted by $w$, since $\alpha$ is a birational morphism. By construction, the restriction of $w$ on $\lfloor B\rfloor$ is equal to $s^{p^{e+e'}}$.

To conclude, for every $s\in \operatorname{A}\bigl(\lfloor B\rfloor, m(K_X+B)|_{\lfloor B\rfloor}\bigl)$, we can find an integer $e+e'$, such that there exists a section $w\in \operatorname{PA}\bigl(X,p^{e+e'}m(K_X+B)\bigl)$, such that
$$s^{p^{e+e'}}=w|_{\lfloor B\rfloor}.$$
This proves the first part of the conclusion. To prove the remaining part, we take sections
$$\{s_i\}_{i=1}^d\subset\operatorname{A}(\lfloor B\rfloor,m(K_X+B)|_{\lfloor B\rfloor})$$
such that $\{s_i\}_{i=1}^d$ generate $\Oo_{\lfloor B\rfloor}(m(K_X+B)|_{\lfloor B\rfloor})$. By the same construction above, we first get
$$\{s_i'':=\alpha|_{B_Y^{=1}}^*(s_i)\}_{i=1}^d\subset H^0(B_Y^{=1},m(K_Y+B_Y)|_{B_Y^{=1}}),$$
which generate $\Oo_{B_Y^{=1}}(m(K_Y+B_Y)|_{B_Y^{=1}})$. And after taking a $p^e$-th power, $\{s_i''^{p^e}\}_{i=1}^d$ descend to sections
$$\{s_i'\}_{i=1}^d\subset H^0(\lfloor B'\rfloor,p^em(K_{X'}+B')|_{\lfloor B'\rfloor})$$
which generate $\Oo_{\lfloor B'\rfloor}(p^em(K_{X'}+B')|_{\lfloor B'\rfloor})$. The restriction
$$\{s_i'|_{B'^h}\}_{i=1}^d\subset\operatorname{A}(B'^h,p^em(K_{X'}+B')|_{B'^h})$$
generate $\Oo_{ B'^h}(p^em(K_{X'}+B')|_{ B'^h})$. And thus
$$\{s_i'^{p^{e'}}|_{B'^h}\}_{i=1}^d\subset\operatorname{A}(B'^h,p^{e+e'}m(K_{X'}+B')|_{B'^h})$$
generate $\Oo_{ B'^h}(p^{e+e'}m(K_{X'}+B')|_{ B'^h})$.

Given the fact that $B'^h$ dominates $V$, we have that the lift of these sections on $X'$
$$\{w'_i\}_{i=1}^d\subset \operatorname{PA}(X',p^{e+e'}m(K_{X'}+B'))$$
generates $\Oo_{X'}(p^{e+e'}m(K_{X'}+B'))$, which implies $\{w_i\}_{i=1}^d$ generates $\Oo_{X}(p^{e+e'}m(K_{X}+B))$.
\end{proof}

The following lemma, corresponding to \cite[Lemma 4.9]{Fuj00}, is used in the proof of Theorem \ref{cor:slc-threefold-abundance} and Theorem \ref{main thm: finiteness of pluri for peff lc surf}  in the case of Calabi-Yau pairs.

\begin{lemma}\label{bir stable SG}
Let $(X,B)$ be a pure $d$-dimensional projective divisorial log terminal pair over an algebraically closed field $k$ of characteristic $p>5$ with $K_X+B$ nef and $d\leq 3$ a positive integer. For a sufficiently large and divisible integer $m$, if $s\in \operatorname{PA}\bigl(X,m(K_X+B)\bigr)$, then
\[
(g_*s)|_{\lfloor B\rfloor}=s|_{\lfloor B\rfloor}
\quad\text{and}\quad
g_*s\in \operatorname{PA}\bigl(X,m(K_X+B)\bigr)
\]
for every $g\in \Bir(X,B)$.
\end{lemma}

\begin{proof}
For a divisor $B=\sum_i b_iB_i$, we write $B^{=1}:=\sum_{b_i=1}B_i$. Let
\[
\alpha,\beta:(Y,B_Y)\longrightarrow (X,B)
\]
be a common log resolution of a $B$-birational map $g:X\dashrightarrow X$ such that $\alpha=g\circ\beta$. To prove $(g_*s)|_{\lfloor B\rfloor}=s|_{\lfloor B\rfloor}$, it is sufficient to check that
\[
\alpha|_{B_Y^{=1}}^*\bigl(s|_{\lfloor B\rfloor}\bigr)
=
\beta|_{B_Y^{=1}}^*\bigl(s|_{\lfloor B\rfloor}\bigr)
\]
in $H^0\!\left(B_Y^{=1},\mathcal{O}_{B_Y^{=1}}\bigl(m(K_Y+B_Y)|_{B_Y^{=1}}\big)\right)$ for some common log resolution $(Y,B_Y)$ since the restricted morphisms $B_Y^{=1}\to\lfloor B\rfloor$ are dominant. For this purpose we first establish the following two claims. For convenience, we write $W$ instead of $(W,\Xi)$ if the pair structure is obvious from the context. Let $n$ be a positive integer at most $3$.

\noindent \textbf{Claim (A\(_n\)).} \textit{Let \((T_1, \Theta_1)\) and \((T_2, \Theta_2)\) be $n$-dimensional log smooth sub-pairs with log canonical singularities and \(h :  (T_1,\Theta_1)\to(T_2,\Theta_2) \) a \(B\)-birational morphism. If \(W_2 \subset T_2\) is a log canonical center of $(T_2, \Theta_2)$, then there is a log canonical center $W_1$ of  $(T_1, \Theta_1)$, such that the restricted morphism \(h|_{W_1}\) is a \(B\)-birational morphism between the pairs on $W_1$ and $W_2$ induced by $(T_1,\Theta_1)$ and $(T_2,\Theta_2)$ respectively from adjunction.}

\vspace{0.5em}

\noindent \textbf{Claim (B\(_n\)).} \textit{Let \((T_1, \Theta_1)\) and \((T_2, \Theta_2)\) be \(n\)-dimensional log smooth sub-pairs with log canonical singularities and \(h : T_1 \to T_2\) a \(B\)-birational morphism. Assume that \(W \subset T_1\) is a log canonical center of $(T_1,\Theta_1)$. If  \(W \to h(W)\) is not \(B\)-birational, then there is a log canonical center \(W' \subseteq W\) of $(T_1,\Theta_1)$, such that \(W' \to h(W)\) is a \(B\)-birational morphism and the inclusion \(W' \to W\) induces an isomorphism \[H^0(W, m(K_{T_1} + \Theta_1)|_W) \simeq H^0(W', m(K_{T_1} + \Theta_1)|_{W'}).\]}

\vspace{0.5em}

\noindent \textit{Proof of Claims.} We prove Claim (A\(_n\)) and Claim (B\(_n\)) by induction on \(n\).  For the case $n=1$, (A\(_1\)) and (B\(_1\)) are immediate.

First, we treat (A\(_n\)). If \(W_2\) is a divisor, then we can take the proper transform of \(W_2\) as \(W_1\). Otherwise, take a divisor \(V_2 \subseteq T_2\) such that \(W _2\subseteq V_2\) and $V_2$ is a log canonical center of $(T_2,\Theta_2)$. We define \(V_1 \subseteq T_1\) as the proper transform of \(V_2\). By using the (A\(_{n-1}\)) to the \(B\)-birational morphism \(V_1 \to V_2\), we obtain \(W_1 \subseteq V_1\) such that \(W_1 \to W_2\) is a \(B\)-birational morphism. Thus, we have (A\(_n\)).

Next, we treat (B\(_n\)) by using (A\(_l\)) with \(l < n\). Let \(u : (T_1', \Theta_1') \to T_1\) be the blowing-up of \(T_1\) with center \(W\). It is enough to prove (B\(_n\)) for \((T_1',\Theta_1')\to (T_2,\Theta_2)\) and the $u$-exceptional divisor \(E \subseteq T_1'\). It is because if (B$_n$), applied on \((T_1',\Theta_1')\to (T_2,\Theta_2)\) and \(E\), gives \(G \subseteq E\), then \(u(G)\) is the required log canonical center. So we may assume that \(W\) is a divisor. By \cite[2.45]{Kollár_Mori_1998}, we have a sequence of blow-ups
\begin{center}
    \begin{tikzcd}
S_0 \arrow[r, "a_1"] \arrow[rr, "a_2", bend left] \arrow[rrr, "a_3", bend right] \arrow[rrrr, "a_k", bend right] & S_1 \arrow[r] & S_2 \arrow[r] & \cdots \arrow[r] & S_k=T_2
\end{tikzcd}
\end{center}
centering at the centers associated to the valuation of \(W\) such that the induced rational map \(f_0: T_1 \dashrightarrow S_0\) is an isomorphism at the generic point of \(W\). We can replace \((T_1, \Theta_1)\) with the elimination of indeterminacy \((\bar T, \bar \Theta)\) of $f_0$, since \(W\) and its proper transform are \(B\)-birational. So we can assume that \(f_0 : T_1 \to S_0\) is a \(B\)-birational morphism. Since we can apply (A\(_{n-1}\)) to the \(B\)-birational morphism \(f_0 : W \to f_0(W)\), it follows that to prove the remaining part of the claim, we only have to prove (B\(_n\)) to the \(B\)-birational morphism \(S_0 \to S_k=T_2\) and $f_0(W)\subset S_0$. This is because if (B$_n$) for $S_0\to S_k$ and $f_0(W)$ gives a $B$-birational morphism $W_0\to h(W)$ with $W_0\subset f_0(W)$, then (A$_{n-1}$) for $W\to f_0(W)$ and $W_0\subset f_0(W)$ gives the required $W'\subset W$. And also, the isomorphism
$$H^0(W, m(K_{T_1} + \Theta_1)|_W) \simeq H^0(W', m(K_{T_1} + \Theta_1)|_{W'})$$
follows from the isomorphisms
$$H^0(W, m(K_{T_1} + \Theta_1)|_W) \simeq H^0(f_0(W), m(K_{S_0} + \Theta_{S_0})|_{f_0(W)}),$$
$$H^0(f_0(W), m(K_{S_0} + \Theta_{S_0})|_{f_0(W)})\simeq H^0(W_0, m(K_{S_0} + \Theta_{S_0})|_{W_0})$$
and
$$H^0(W_0, m(K_{S_0} + \Theta_{S_0})|_{W_0})\simeq H^0(W', m(K_{T_1} + \Theta_{1})|_{W'}),$$
where $(S_i,\Theta_{S_i})$ are the induced pairs on $S_i$ for $i=0,1,\cdots,k$.

We take a divisor \(D_1 \subseteq S_1\) such that $D_1$ is a log canonical center and \(h(W) \subseteq D_1\). Denote the proper transform of $D_1$ to be \(D_0 \subseteq S_0\). By applying (B\(_{n-1}\)) to $D_0\to D_1$ and \(D_0 \cap f_0(W) \subseteq D_0\), we find a log canonical center $W_{0,1}\subset D_0 \cap f_0(W)\subset f_0(W)$ with a $B$-birational morphism $W_{0,1}\to a_1(f_0(W))$. The isomorphism
$$H^0(f_0(W), m(K_{S_0} + \Theta_{S_0})|_{f_0(W)})\simeq H^0(W_{0,1}, m(K_{S_0} + \Theta_{S_0})|_{W_{0,1}})$$
is given by
$$H^0(f_0(W), m(K_{S_0} + \Theta_{S_0})|_{f_0(W)})\simeq H^0(a_1(f_0(W)), m(K_{S_1} + \Theta_{S_1})|_{a_1(f_0(W))})$$
since $a_1$ is a blow-up centering at $a_1(f_0(W))$ with $f_0(W)$ the exceptional divisor, and
$$H^0(a_1(f_0(W)), m(K_{S_1} + \Theta_{S_1})|_{a_1(f_0(W))})\simeq H^0(W_{0,1}, m(K_{S_0} + \Theta_{S_0})|_{W_{0,1}}).$$
Similarly we obtain a log canonical center $W_{1,1}\subset a_1(f_0(W))$ with a $B$-birational morphism $W_{1,1}\to a_2(f_0(W))$. Applying (A$_{l}$) to the $B$-birational morphism $W_{0,1}\to a_1(f_0(W))$ with $W_{1,1}\subset a_1(f_0(W))$, we obtain a log canonical center $W_{0,2}\subset W_{0,1}\subset f_0(W)$ with a $B$-birational morphism $W_{0,2}\to W_{1,1}$. Thus, we know $W_{0,2}\to a_2(f_0(W))$ is a $B$-birational morphism. Continuing the process, we can find a log canonical center $W_{0,k}\subset f_0(W)$ with $W_{0,k}\to a_k(f_0(W))=h(W)$ $B$-birational. And $W_{0,k}$ is the required $W_0$. The isomorphism
$$H^0(f_0(W), m(K_{S_0} + \Theta_{S_0})|_{f_0(W)})\simeq H^0(W_0, m(K_{S_0} + \Theta_{S_0})|_{W_0})$$
follows by construction.\hfill$\square$
\vspace{0.5em}

\noindent \textit{Proof of Lemma \ref{bir stable SG} continued.} Let \(f : (X', B') \to (X, B)\) be a log resolution such that \(K_{X'} + B' = f^*(K_X + B)\) with the coefficient $\text{Coeff}_{B'}(E)<1$ for any $f$-exceptional divisor $E$ (see Definition \ref{sing of pairs}). Let \(\alpha = \amalg \alpha_i\) and \(\beta = \amalg \beta_i : (Y, \Theta) = \amalg (Y_i, \Theta_i) \to (X, B)\) be a common log resolution of a \(B\)-birational map \(g : X\dashrightarrow X\) passing through \((X', B')\) such that \(\alpha = g \circ \beta\). We write \(\alpha', \beta' : Y \to X'\) such that \(\alpha = f \circ \alpha', \beta = f \circ \beta'\).

For any irreducible divisor \(D \subset \Theta^{=1}\). Apply (B\(_n\)) to the \(\alpha'_i : Y_i \to X'\) with \(D \subseteq Y_i\). Then we obtain a log canonical center \(D' \subseteq D\) or \(D' = D\) such that \(D' \to \alpha'_i(D)\) is \(B\)-birational. Apply (B\(_n\)) to the \(\beta'_i : Y_i \to X'\) with $D'\subset Y_i$. Then we obtain a log canonical center \(D'' \subseteq D'\) or \(D'' = D'\) such that \(D'' \to \beta'_i(D')\) is \(B\)-birational. By repeating the above construction, we obtain \(F \subseteq D\) such that \(\alpha' : F \to \alpha'(F)\) and \(\beta' : F \to \beta'(F)\) are \(B\)-birational. The morphism \(f : \alpha'(F) \to \alpha(F)\) and \(f : \beta'(F) \to \beta(F)\) are \(B\)-birational since \(f\) induces an isomorphism on the generic point of any log canonical center of $(X',B')$. Then \(\alpha : F \to \alpha(F)\) and \(\beta : F \to \beta(F)\) are \(B\)-birational. Since \(s \in \mathrm{PA}(X, m(K_X + B))\), it follows that \(\alpha^*(s|_{\alpha(F)}) = \beta^*(s|_{\beta(F)})\). Given \(H^0(D, m(K_Y + \Theta)|_D) \simeq H^0(F, m(K_Y + \Theta)|_F)\) implied by (B\(_n\)), we have \(\alpha^*s = \beta^*s\) on \(D\). Since \(D\) is an arbitrary irreducible component of \(\Theta^{=1}\), we obtain \(\alpha^*(s|_{\llcorner B \lrcorner}) = \beta^*(s|_{\llcorner B \lrcorner})\). Thus \(g_*s|_{\llcorner B \lrcorner} = s|_{\llcorner B \lrcorner}\).
\end{proof}

Let us prove the main result of this paper.

\begin{thm}\label{cor:slc-threefold-abundance}
Let $(X, B)$ be a projective semi log canonical pair of dimension $3$ over an algebraically closed field $k$ of characteristic $p>42$. If $K_X + B$ is nef, then it is semi-ample.
\end{thm}

\begin{proof}
We first prove the following claims.

\noindent \textbf{Claim (A\(_n\)).} Let $(X, B)$ be a projective $\Qq$-factorial divisorial log terminal pair of dimension $n\le 3$ over an algebraically closed field $k$ of characteristic $p>42$. If $K_X+B$ is nef, then $\operatorname{PA}(X, m(K_X+B))$ generates $\Oo_X(m(K_X+B))$ for a sufficiently divisible positive integer $m$.

\noindent \textbf{Claim (B\(_n\)).} Let $(X,B)$ be a projective semi divisorial log terminal pair of dimension $n\le 2$ over an algebraically closed field $k$ of characteristic $p>42$. If $K_X+B$ is nef, then $\operatorname{A}(X, m(K_X+B))$ generates $\Oo_X(m(K_X+B))$ for a sufficiently divisible positive integer $m$.

\noindent \textit{Proof of Claims.} We prove these claims inductively. Claim $(B_0)$ is immediate.
\begin{itemize}
    \item Claim $(B_{n-1})$ implies Claim $(A_n)$ for $n\le 3$:
\end{itemize}

This follows from Proposition \ref{prop:preadmissible-extension}. Note that $(\lfloor B\rfloor,\Diff(B-\lfloor B\rfloor))$ is semi divisorial log terminal (see Remark \ref{sdlt of reduced part of dlt 3-folds}).

\begin{itemize}
\item      Claim $(A_n)$ implies Claim $(B_n)$ for $n\le 2$:
\end{itemize}

Take the normalization $X':=\coprod_iX_i\to X=\cup_i X_i$. Let $(X',B')=\coprod_i(X_i,B_i')$ be the pair defined as in Definition \ref{def : sdlt slc}. Then $(X',B')$ is divisorial log terminal. Define a group $G:=\rho_{m}(\Bir(X',B'))$ to be the image of the $B$-pluricanonical representation of $(X',B')$ for a sufficiently divisible positive integer $m$. By Theorem \ref{1-dim subklt fin rep} and Theorem \ref{main thm: finiteness of pluri for peff lc surf & Cartier ind}, $G$ is finite. Write $G=\{g_j\}_{j=1}^N$.

Take sections $\{u_i\}_{i=1}^M\subset \operatorname{PA}(X',m(K_{X'}+B'))$ such that $\bigcap_{i}\{u_i=0\}$ is empty. By Lemma \ref{bir stable SG}, we have $g_ju_i|_{\lfloor B'\rfloor}=u_i|_{\lfloor B'\rfloor}$ and thus $g_ju_i\in \operatorname{PA}(X',m(K_{X'}+B'))$ for all $1\le j\le N$ and $1\le i\le M$. Therefore, we know that
\[
s_{ij}:=\sigma_j^{N!/j}(g_1u_i,\ldots,g_Nu_i)
\in
\operatorname{A}(X',mN!(K_{X'}+B')),
\]
where $\sigma_j$ is the $N$-variable elementary symmetric polynomial of degree $j$. Then by the following fact
\[
\{u_i=0\}
\supset
\{g_1^*u_i=\ldots=g_N^*u_i=0\}
=
\bigcap_{j=1}^{N}
\{\sigma_j(g_1^*u_i,\ldots,g_N^*u_i)=0\},
\]
we obtain that $\bigcap_{i,j}\{s_{ij}=0\}$ is empty, and thus we can prove that $\operatorname{A}(X',mN!(K_{X'}+B'))$ generates $\mathcal{O}_{X'}(mN!(K_{X'}+B'))$. By Lemma \ref{lem:pa-descends-normalization}, and the definition of admissible sections (Definition \ref{def:admissible-section}), we know that $\{\nu_*s_{ij}\}_{i,j}\subset \operatorname{A}(X,mN!(K_X+B))$ and generates $\mathcal{O}_X(mN!(K_{X}+B))$.\hfill$\square$

\noindent \textit{Proof of Theorem \ref{cor:slc-threefold-abundance} continued.} We take the normalization of $(X,B)$ and take a $\mathbb{Q}$-factorial divisorial log terminal model $(Y,\Xi)=\coprod (Y_i,\Xi_i)$ for each irreducible component by \cite[Theorem 1.6]{birkar2016existence}. By Lemma \ref{lem:pa-descends-normalization}, it remains to prove $\operatorname{PA}(Y,m(K_Y+\Xi))$ generates $\mathcal{O}_Y(m(K_Y+\Xi))$ for a sufficiently divisible integer $m>0$. Thus we can conclude by applying Claim $(A_3)$.
\end{proof}

\bibliographystyle{alpha}
\bibliography{refs}

@Inbook{Ueno1975,
author="Ueno, Kenji",
title="Addition formula for Kodaira dimensions of analytic fibre bundles",
bookTitle="Classification Theory of Algebraic Varieties and Compact Complex Spaces: Notes written in collaboration with P. Cherenack",
year="1975",
publisher="Springer Berlin Heidelberg",
address="Berlin, Heidelberg",
pages="172--187",
isbn="978-3-540-37415-2",
doi="10.1007/BFb0070576",
url="https://doi.org/10.1007/BFb0070576"
}

@article{Fuj00,
  author = {Osamu Fujino},
  title = {Abundance theorem for semi log canonical threefolds},
  journal = {Duke Mathematical Journal},
  volume = {102},
  number = {3},
  pages = {513--532},
  year = {2000},
  url = {https://www.math.kyoto-u.ac.jp/~fujino/Abundance.pdf},
}

@article{NT20,
  author  = {Nakamura, Yusuke and Tanaka, Hiromu},
  title   = {A {Witt} {Nadel} vanishing theorem for threefolds},
  journal = {Compositio Mathematica},
  volume  = {156},
  number  = {3},
  pages   = {435--475},
  year    = {2020},
  doi     = {10.1112/S0010437X1900770X}
}

@misc{Gu2020SurfaceFW,
  author = {Gu, Yi},
  title = {Surface fibrations with large equivariant automorphism group},
  year = {2020},
  howpublished = {arXiv preprint arXiv:2004.12656},
  eprint = {2004.12656},
  archivePrefix = {arXiv},
  primaryClass = {math.AG},
  doi = {10.48550/arXiv.2004.12656},
  url = {https://arxiv.org/abs/2004.12656},
}

@article{Chen_2024,
  author = {Yifei Chen and Constantin Shramov},
  title = {Automorphisms of surfaces over fields of positive characteristic},
  journal = {Geometry \& Topology},
  volume = {28},
  number = {6},
  pages = {2747--2791},
  year = {2024},
  doi = {10.2140/gt.2024.28.2747},
  url = {https://msp.org/gt/2024/28-6/gt-v28-n6-p04-p.pdf},
}

@misc{Pat17,
  author        = {Patakfalvi, Zsolt},
  title         = {On the projectivity of the moduli space of stable surfaces
                   in characteristic {$p>5$}},
  year          = {2017},
  eprint        = {1710.03818v1},
  archivePrefix = {arXiv},
  primaryClass  = {math.AG},
  url           = {https://arxiv.org/abs/1710.03818v1}
}

@incollection{hacon2012finitenessbrepresentationsemilogcanonical,
  author = {Christopher D. Hacon and Chenyang Xu},
  editor = {J{\'a}nos Koll{\'a}r and Osamu Fujino and Shigeru Mukai and Noboru Nakayama},
  title = {On finiteness of {B}-representations and semi-log canonical abundance},
  booktitle = {Minimal Models and Extremal Rays (Kyoto, 2011)},
  series = {Advanced Studies in Pure Mathematics},
  volume = {70},
  pages = {361--377},
  publisher = {Mathematical Society of Japan},
  address = {Tokyo},
  year = {2016},
  doi = {10.2969/aspm/07010361},
  url = {https://doi.org/10.2969/aspm/07010361},
}

@article{Fujino_Gongyo_2014,
  author = {Fujino, Osamu and Gongyo, Yoshinori},
  title = {Log pluricanonical representations and the abundance conjecture},
  journal = {Compositio Mathematica},
  volume = {150},
  number = {4},
  pages = {593--620},
  year = {2014},
  doi = {10.1112/S0010437X13007495},
  url = {https://www.cambridge.org/core/journals/compositio-mathematica/article/abs/log-pluricanonical-representations-and-the-abundance-conjecture/DF078000E896B4D757027B284ADF8DC1},
}

@incollection{kollár2012sourceslogcanonicalcenters,
  author = {J{\'a}nos Koll{\'a}r},
  editor = {J{\'a}nos Koll{\'a}r and Osamu Fujino and Shigeru Mukai and Noboru Nakayama},
  title = {Sources of log canonical centers},
  booktitle = {Minimal Models and Extremal Rays (Kyoto, 2011)},
  series = {Advanced Studies in Pure Mathematics},
  volume = {70},
  pages = {29--48},
  publisher = {Mathematical Society of Japan},
  address = {Tokyo},
  year = {2016},
  url = {https://arxiv.org/abs/1107.2863},
}

@book{kollar-moduli-hd,
  author = {J{\'a}nos Koll{\'a}r},
  title = {Families of Varieties of General Type},
  series = {Cambridge Tracts in Mathematics},
  volume = {231},
  publisher = {Cambridge University Press},
  year = {2023},
  doi = {10.1017/9781009346115},
  url = {https://www.cambridge.org/core/books/families-of-varieties-of-general-type/423EBD2D90435E5EA8AAC69CE0E6EFD1},
  note = {With the collaboration of Klaus Altmann and S{\'a}ndor J. Kov{\'a}cs},
}

@article{Jiang_2021,
  author = {Chen Jiang and Haidong Liu},
  title = {Boundedness of log pluricanonical representations of log {Calabi--Yau} pairs in dimension 2},
  journal = {Algebra \& Number Theory},
  volume = {15},
  number = {2},
  pages = {545--567},
  year = {2021},
  doi = {10.2140/ant.2021.15.547},
  url = {https://msp.org/ant/2021/15-2/p07.xhtml},
}

@article{Nygaard1980HigherDR,
  title={Higher de Rham-Witt complexes of supersingular K3 surfaces},
  author={Niels O. Nygaard},
  journal={Compositio Mathematica},
  year={1980},
  volume={42},
  pages={245-271},
  url={https://api.semanticscholar.org/CorpusID:54986113}
}

@article{jang2016representationsautomorphismgroupsfrobenius,
  author = {Junmyeong Jang},
  title = {The representations of the automorphism groups and the {Frobenius} invariants of {K3} surfaces},
  journal = {Michigan Mathematical Journal},
  volume = {65},
  number = {1},
  pages = {147--163},
  year = {2016},
  doi = {10.1307/mmj/1457101815},
  url = {https://doi.org/10.1307/mmj/1457101815},
}

@book{Kollár_Mori_1998,
  author = {J{\'a}nos Koll{\'a}r and Shigefumi Mori},
  title = {Birational Geometry of Algebraic Varieties},
  series = {Cambridge Tracts in Mathematics},
  volume = {134},
  publisher = {Cambridge University Press},
  year = {1998},
  doi = {10.1017/CBO9780511662560},
  url = {https://www.cambridge.org/core/books/birational-geometry-of-algebraic-varieties/FB73FD7FDC65F023521ED6B1BB7F3EC7},
}

@article{Matsumoto_2023,
  author = {Yuya Matsumoto},
  title = {Degeneration of {K3} surfaces with non-symplectic automorphisms},
  journal = {Rendiconti del Seminario Matematico della Universit{\`a} di Padova},
  volume = {150},
  pages = {227--245},
  year = {2023},
  doi = {10.4171/RSMUP/123},
  url = {https://ems.press/content/serial-article-files/44136},
}

@article{tanaka2015abundancetheoremsemilog,
  author  = {Tanaka, Hiromu},
  title   = {Abundance theorem for semi log canonical surfaces in positive characteristic},
  journal = {Osaka Journal of Mathematics},
  year    = {2016},
  volume  = {53},
  number  = {2},
  pages   = {535--566},
  doi     = {10.18910/58879}
}

@book{Kollár_2013,
  author = {J{\'a}nos Koll{\'a}r},
  title = {Singularities of the Minimal Model Program},
  series = {Cambridge Tracts in Mathematics},
  volume = {200},
  publisher = {Cambridge University Press},
  year = {2013},
  doi = {10.1017/CBO9781139547895},
  url = {https://assets.cambridge.org/97811070/35348/frontmatter/9781107035348_frontmatter.pdf},
  note = {With a collaboration of S{\'a}ndor Kov{\'a}cs},
}

@incollection{Milne1986,
  author = {J. S. Milne},
  editor = {Gary Cornell and Joseph H. Silverman},
  title = {Abelian Varieties},
  booktitle = {Arithmetic Geometry},
  pages = {103--150},
  publisher = {Springer},
  year = {1986},
  doi = {10.1007/978-1-4613-8655-1_5},
  url = {https://link.springer.com/chapter/10.1007/978-1-4613-8655-1_5},
}

@article{abundanceforlcsurf,
author = {Hiromu Tanaka},
title = {{Minimal models and abundance for positive characteristic log surfaces}},
volume = {216},
journal = {Nagoya Mathematical Journal},
number = {none},
publisher = {Nagoya Mathematical Journal},
pages = {1 -- 70},
year = {2014},
doi = {10.1215/00277630-2801646},
URL = {https://doi.org/10.1215/00277630-2801646}
}

@article{Nakamura_Ueno[1],
author = {Iku NAKAMURA and Kenji UENO},
title = {{An addition formula for Kodaira dimensions of analytic fibre bundles whose fibres are Moi\v{s}ezon manifolds}},
volume = {25},
journal = {Journal of the Mathematical Society of Japan},
number = {3},
publisher = {Mathematical Society of Japan},
pages = {363 -- 371},
year = {1973},
doi = {10.2969/jmsj/02530363},
URL = {https://doi.org/10.2969/jmsj/02530363}
}

@article{Posva_2023_slc,
   title={Abundance for slc surfaces over arbitrary fields},
   volume={Volume 7},
   ISSN={2491-6765},
   url={http://dx.doi.org/10.46298/epiga.2023.volume7.8803},
   DOI={10.46298/epiga.2023.volume7.8803},
   journal={Épijournal de Géométrie Algébrique},
   publisher={Centre pour la Communication Scientifique Directe (CCSD)},
   author={Posva, Quentin},
   year={2023},
   month=Feb }

@article{Filipazzi2020ConnectednessPF,
  author = {Filipazzi, Stefano and Waldron, Joe},
  title = {Connectedness principle for 3-folds in characteristic $p>5$},
  journal = {Michigan Mathematical Journal},
  volume = {74},
  number = {4},
  pages = {675--701},
  year = {2024},
  doi = {10.1307/mmj/20216143},
  url = {https://arxiv.org/abs/2010.08414},
}

@article{katsura1984elliptic,
  author = {Toshiyuki Katsura and Kenji Ueno},
  title = {On elliptic surfaces in characteristic $p$},
  journal = {Mathematische Annalen},
  volume = {272},
  number = {3},
  pages = {291--330},
  year = {1985},
  doi = {10.1007/BF01455561},
  url = {https://link.springer.com/article/10.1007/BF01455561},
}

@book{buadescu2001algebraic,
  author = {B{\u a}descu, Lucian},
  title = {Algebraic Surfaces},
  series = {Universitext},
  publisher = {Springer-Verlag},
  address = {New York},
  year = {2001},
  doi = {10.1007/978-1-4757-3512-3},
  url = {https://link.springer.com/book/10.1007/978-1-4757-3512-3},
}

@article{bejleri2026moduli,
  author = {Bejleri, Dori and Foster, Josiah and Fernandez Herrero, Andres and Inchiostro, Giovanni and Makarova, Svetlana and Zhao, Junyan},
  title = {Moduli of elliptic surfaces of {Kodaira} dimension one fibered over rational curves},
  journal = {Advances in Mathematics},
  volume = {494},
  pages = {110942},
  year = {2026},
  doi = {10.1016/j.aim.2026.110942},
  url = {https://www.sciencedirect.com/science/article/pii/S0001870826001647},
}

@incollection{liedtke2013algebraic,
  author = {Christian Liedtke},
  title = {Algebraic Surfaces in Positive Characteristic},
  booktitle = {Birational Geometry, Rational Curves, and Arithmetic},
  series = {Simons Symposia},
  pages = {229--292},
  publisher = {Springer},
  year = {2013},
  doi = {10.1007/978-1-4614-6482-2_11},
  url = {https://link.springer.com/chapter/10.1007/978-1-4614-6482-2_11},
}

@book{milne2017algebraicgroups,
  author = {J. S. Milne},
  title = {Algebraic Groups: The Theory of Group Schemes of Finite Type over a Field},
  series = {Cambridge Studies in Advanced Mathematics},
  volume = {170},
  publisher = {Cambridge University Press},
  year = {2017},
  doi = {10.1017/9781316711736},
  url = {https://www.cambridge.org/core/books/algebraic-groups/48D122BAEC63C2B70304E20F34482C8F},
}

@incollection{Bombieri_Mumford_1977,
  author = {Enrico Bombieri and David Mumford},
  editor = {Baily, Jr., Walter L. and Shioda, Tetsuji},
  title = {{Enriques}' Classification of Surfaces in Char. $p$, {II}},
  booktitle = {Complex Analysis and Algebraic Geometry: A Collection of Papers Dedicated to K. {Kodaira}},
  pages = {23--42},
  publisher = {Cambridge University Press},
  year = {1977},
  doi = {10.1017/CBO9780511569197.004},
  url = {https://www.cambridge.org/core/books/abs/complex-analysis-and-algebraic-geometry/enriques-classification-of-surfaces-in-char-p-ii/6B32C9CE0093D341878EA19D2C030B78},
}

@article{PW17,
  title={Singularities of general fibers and the LMMP},
  author={Zsolt Patakfalvi and Joe Waldron},
  journal={American Journal of Mathematics},
  year={2017},
  volume={144},
  pages={505 - 540},
  url={https://api.semanticscholar.org/CorpusID:119158578}
}

@article{EGAIV,
  author  = {Grothendieck, Alexander},
  title   = {{\'E}l{\'e}ments de g{\'e}om{\'e}trie alg{\'e}brique.
             {IV}. {\'E}tude locale des sch{\'e}mas et des morphismes
             de sch{\'e}mas},
  journal = {Publications Math{\'e}matiques de l'{IH{\'E}S}},
  volume  = {20, 24, 28, 32},
  year    = {1964--1967},
  note    = {Parts I--IV}
}

@article{haconxu2015,
  title={On the three dimensional minimal model program in positive characteristic},
  author={Hacon, Christopher and Xu, Chenyang},
  journal={Journal of the American Mathematical Society},
  volume={28},
  number={3},
  pages={711--744},
  year={2015}
}

@article{cascini2015base,
  author = {Cascini, Paolo and Tanaka, Hiromu and Xu, Chenyang},
  title = {On base point freeness in positive characteristic},
  journal = {Annales scientifiques de l'{\'E}cole Normale Sup{\'e}rieure},
  volume = {48},
  number = {5},
  pages = {1239--1272},
  year = {2015},
  doi = {10.24033/asens.2269},
  url = {https://smf.emath.fr/publications/sur-la-vacuite-du-lieu-base-en-caracteristique-positive},
}

@article{birkar2016existence,
  author = {Birkar, Caucher},
  title = {Existence of flips and minimal models for 3-folds in char. $p$},
  journal = {Annales scientifiques de l'{\'E}cole Normale Sup{\'e}rieure},
  volume = {49},
  number = {1},
  pages = {169--212},
  year = {2016},
  doi = {10.24033/asens.2279},
  url = {https://www.numdam.org/articles/10.24033/asens.2279/},
}

@article{bernasconi2023vanishing,
  author = {Bernasconi, Fabio and Koll{\'a}r, J{\'a}nos},
  title = {Vanishing theorems for three-folds in characteristic $p>5$},
  journal = {International Mathematics Research Notices},
  volume = {2023},
  number = {4},
  pages = {2846--2866},
  year = {2023},
  doi = {10.1093/imrn/rnab316},
  url = {https://academic.oup.com/imrn/article-abstract/2023/4/2846/6433424},
}

@article{Brion2021,
  author = {Michel Brion},
  title = {Some automorphism groups are linear algebraic},
  journal = {Moscow Mathematical Journal},
  volume = {21},
  number = {3},
  pages = {453--466},
  year = {2021},
  doi = {10.17323/1609-4514-2021-21-3-453-466},
  url = {https://www.mathnet.ru/php/archive.phtml?jrnid=mmj&option_lang=eng&paperid=801&wshow=paper},
}

@misc{bernasconi2026abundancethreefoldspositivecharacteristic,
  author = {Fabio Bernasconi and Jakub Witaszek and Zheng Xu},
  title = {Abundance for threefolds in positive characteristic},
  year = {2026},
  howpublished = {arXiv preprint arXiv:2610.10439},
  eprint = {2610.10439},
  archivePrefix = {arXiv},
  primaryClass = {math.AG},
  doi = {10.48550/arXiv.2610.10439},
  url = {https://arxiv.org/abs/2610.10439},
}

@article{kollar1988threefolds,
  author = {J{\'a}nos Koll{\'a}r and Nicholas I. Shepherd-Barron},
  title = {Threefolds and deformations of surface singularities},
  journal = {Inventiones Mathematicae},
  volume = {91},
  number = {2},
  pages = {299--338},
  year = {1988},
  doi = {10.1007/BF01389370},
  url = {https://link.springer.com/article/10.1007/BF01389370},
}

@article{abundanceforthreefolds,
author = {Sean Keel and Kenji Matsuki and James McKernan},
title = {{Log abundance theorem for threefolds}},
volume = {75},
journal = {Duke Mathematical Journal},
number = {1},
publisher = {Duke University Press},
pages = {99 -- 119},
year = {1994},
doi = {10.1215/S0012-7094-94-07504-2},
URL = {https://doi.org/10.1215/S0012-7094-94-07504-2}
}

@book{Kollr1992FlipsAA,
  editor = {J{\'a}nos Koll{\'a}r},
  title = {Flips and abundance for algebraic threefolds},
  series = {Ast{\'e}risque},
  volume = {211},
  publisher = {Soci{\'e}t{\'e} Math{\'e}matique de France},
  address = {Paris},
  year = {1992},
  doi = {10.24033/ast.198},
  url = {https://numdam.org/issues/AST_1992__211__1_0/},
}

@article{hacon2013birational,
  author = {Christopher D. Hacon and James McKernan and Chenyang Xu},
  title = {On the birational automorphisms of varieties of general type},
  journal = {Annals of Mathematics},
  volume = {177},
  number = {3},
  pages = {1077--1111},
  year = {2013},
  doi = {10.4007/annals.2013.177.3.6},
  url = {https://annals.math.princeton.edu/2013/177-3/p06},
}

@article{hacon2018boundedness,
  author = {Christopher D. Hacon and James McKernan and Chenyang Xu},
  title = {Boundedness of moduli of varieties of general type},
  journal = {Journal of the European Mathematical Society},
  volume = {20},
  number = {4},
  pages = {865--901},
  year = {2018},
  doi = {10.4171/JEMS/778},
  url = {https://ems.press/journals/jems/articles/15330},
}

@book{DemazureGabriel1970,
  author = {Michel Demazure and Pierre Gabriel},
  title = {Groupes alg{\'e}briques. Tome {I}: G{\'e}om{\'e}trie alg{\'e}brique, g{\'e}n{\'e}ralit{\'e}s, groupes commutatifs},
  publisher = {Masson \& Cie; North-Holland Publishing Co.},
  address = {Paris; Amsterdam},
  year = {1970},
}

@book{SGA3,
  editor = {Michel Demazure and Alexander Grothendieck},
  title = {Sch{\'e}mas en groupes {I}, {II}, {III} ({SGA} 3)},
  series = {Lecture Notes in Mathematics},
  volume = {151--153},
  publisher = {Springer-Verlag},
  address = {Berlin--Heidelberg--New York},
  year = {1970},
  url = {https://smf.emath.fr/publications/schemas-en-groupes-sga-3-tome-i-proprietes-generales-des-schemas-en-groupes-seminaire},
}

@article{PS09,
  author  = {Prokhorov, Yuri G. and Shokurov, Vyacheslav V.},
  title   = {Towards the second main theorem on complements},
  journal = {Journal of Algebraic Geometry},
  volume  = {18},
  number  = {1},
  pages   = {151--199},
  year    = {2009},
  doi     = {10.1090/S1056-3911-08-00498-0}
}

@article{Witaszek2021,
  author        = {Witaszek, Jakub},
  title         = {On the canonical bundle formula and log abundance
                   in positive characteristic},
  journal       = {Mathematische Annalen},
  volume        = {381},
  number        = {3-4},
  pages         = {1309--1344},
  year          = {2021},
  doi           = {10.1007/s00208-021-02231-5},
  eprint        = {1711.04380},
  archivePrefix = {arXiv},
  primaryClass  = {math.AG}
}

@article{DH16,
  author = {Das, Omprokash and Hacon, Christopher D.},
  title = {On the adjunction formula for 3-folds in characteristic $p>5$},
  journal = {Mathematische Zeitschrift},
  volume = {284},
  number = {1--2},
  pages = {255--269},
  year = {2016},
  doi = {10.1007/s00209-016-1655-4},
  url = {https://link.springer.com/article/10.1007/s00209-016-1655-4},
}

@unpublished{BWX26indices,
  author = {Bernasconi, Fabio and Witaszek, Jakub and Xu, Zheng},
  title  = {Bounds on the index of log {Calabi--Yau} surfaces},
  note   = {In preparation},
  year   = {2026}
}

\end{document}